\documentclass[a4paper,oneside,11pt]{article}
\usepackage[a4paper]{geometry}
\usepackage{aeguill}                 
\usepackage{graphicx}
\usepackage{caption}
\usepackage{amsmath,amsfonts,amssymb,amsthm}
\usepackage{mathabx}
\usepackage{pstricks,comment} 
\usepackage{pst-plot}
\usepackage{pstricks-add}
\usepackage{multido}
\usepackage{url}
\usepackage{epstopdf,enumerate}
\usepackage{srcltx}
\usepackage[utf8]{inputenc} 
\usepackage[T1]{fontenc} 
\usepackage[english]{babel} 
\usepackage{hyperref}
\usepackage[all]{xy} 
\usepackage{titlesec}
\usepackage{mathabx}
\usepackage{tikz}
\usepackage{fancyhdr}
\usepackage{mathrsfs}
\usepackage[us,12hr]{datetime}
\fancypagestyle{plain}{
\fancyhf{}
\rfoot{ {\ddmmyyyydate\today}}
\lfoot{\thepage}
}

    \newtheoremstyle{TheoremNum}
        {\topsep}{\topsep}              
        {\itshape}                      
        {}                              
        {\bfseries}                     
        {.}                             
        { }                             
        {\thmname{#1}\thmnote{ \bfseries #3}}

{\theoremstyle {definition} \newtheorem {defi} {Definition}[section]}
{\theoremstyle {plain}  \newtheorem {theo} [defi] {Theorem}}
{\theoremstyle {plain}  \newtheorem {coro} [defi] {Corollary}}
{\theoremstyle {plain} \newtheorem {prop} [defi] {Proposition}}
{\theoremstyle {plain} \newtheorem {lem}[defi] {Lemma}}
{\theoremstyle {plain} }
{\theoremstyle {plain} }
{\theoremstyle{TheoremNum} }
{\theoremstyle{TheoremNum} }
{\theoremstyle{TheoremNum} }

\newcommand{\n}{{\tt n}}

\newcommand{\Aut}{\mathrm{Aut}}
\newcommand{\Out}{\mathrm{Out}}

\newcommand{\Mod}{\mathrm{Mod}}

\newcommand{\Stab}{\mathrm{Stab}}
\newcommand{\Poly}{\mathrm{Poly}}

\newcommand{\Curr}{\mathrm{Curr}}
\newcommand{\PCurr}{\mathbb{P}\mathrm{Curr}}
\newcommand{\Supp}{\mathrm{Supp}}
\newcommand{\IA}{\mathrm{IA}}
\newcommand{\id}{\mathrm{id}}

\newcommand{\ZZ}{\mathbb{Z}}

\newcommand{\NN}{\mathbb{N}}
\newcommand{\RR}{\mathbb{R}}

\newcommand{\dem}{\noindent{\bf Proof. }}

\title{Polynomial growth and subgroups of $\Out(F_{\tt n})$}
\author{Yassine Guerch}
\date{\today}
\begin{document}
\maketitle
\renewcommand*\labelenumi{(\theenumi)}

\begin{abstract}
This paper, which is the last of a series of three papers, studies dynamical properties of elements of $\Out(F_{\tt n})$, the outer automorphism group of a nonabelian free group $F_{\tt n}$. We prove that, for every subgroup $H$ of $\Out(F_{\n})$, there exists an element $\phi \in H$ such that, for every element $g$ of $F_{\n}$, the conjugacy class $[g]$ has polynomial growth under iteration of $\phi$ if and only if $[g]$ has polynomial growth under iteration of every element of $H$. 
\end{abstract}

\section{Introduction}\label{Section Introduction}

Let $\n \geq 3$. This paper, which is the last of a series of three papers, studies the exponential growth of elements in $\Out(F_{\n})$. An outer automorphism $\phi \in \Out(F_{\tt n})$ is \emph{exponentially growing} if there exist $g \in F_{\tt n}$, a representative $\Phi$ of $\phi$, a free basis $\mathfrak{B}$ of $F_{\tt n}$ and a constant $K>0$ such that, for every $m \in \NN^*$, we have \begin{equation}\label{Equation intro}
\ell_{\mathfrak{B}}(\Phi^m(g)) \geq e^{Km},
\end{equation}
 where $\ell_{\mathfrak{B}}(\Phi^m(g))$ denotes the length of $\Phi^m(g)$ in the basis $\mathfrak{B}$. If $g \in F_{\n}$ satisfies Equation~\eqref{Equation intro} for every representative $\Phi$ of $\phi$, then $g$ is said to be \emph{exponentially growing under iteration of $\phi$}. Otherwise, one can show, using for instance the technology of relative train tracks introduced by Bestvina and Handel~\cite{BesHan92}, that $g$ has \emph{polynomial growth under iteration of $\phi$}, replacing $\geq e^{Km}$ by $\leq (m+1)^K$ in Equation~\eqref{Equation intro} (see also \cite{Levitt09} for a complete description of all growth types that can occur under iteration of an outer automorphism $\phi$). We denote by $\Poly(\phi)$ the set of elements of $F_{\n}$ which have polynomial growth under iteration of $\phi$. If $H$ is a subgroup of $F_{\n}$, we set $\Poly(H)=\bigcap_{\phi \in H} \Poly(\phi)$. In this article, we prove the following theorem.

\begin{theo}\label{Theo intro 1}
Let ${\tt n} \geq 3$ and let $H$ be a subgroup of $\Out(F_{\tt n})$. There exists $\phi \in H$ such that $\mathrm{Poly}(\phi)=\mathrm{Poly}(H)$.
\end{theo}

In other words, there exists an element of $H$ which encaptures all the exponential growth of $H$: there exists $\phi \in H$ such that if $g \in F_{\n}$ has exponential growth for some element of $H$, then $g$ has exponential growth for $\phi$. Theorem~\ref{Theo intro 1} has analogues in other contexts. For instance, one has a similar result in the context of the mapping class group of a closed, connected, orientable surface $S$ equipped with a hyperbolic structure. Indeed, a consequence of the Nielsen-Thurston classification (see for instance~\cite[Theorem~13.2]{FarMar12}) and the work of Thurston~\cite[Proposition~9.21]{FatLauPoe79} is that the growth of the length of the geodesic representative of a homotopy class of an essential closed curve under iteration of an element of $\Mod(S)$ is either exponential or linear. Moreover, linear growth comes from twists about essential curves while exponential growth comes from pseudo-Anosov homeomorphisms of subsurfaces of $S$. In~\cite{Ivanov92} (see also the work of McCarthy~\cite{McCarthy85}), Ivanov proved that, for every subgroup $H$ of $\Mod(S)$, up to taking a finite index subgroup of $H$, there exists finitely many homotopy classes of pairwise disjoint essential closed curves $C_1,\ldots,C_k$ elementwise fixed by $H$ and such that, for every connected component $S'$ of $S-\bigcup_{i=1}^k C_i$, the restriction \mbox{$H|_{S'} \subseteq \Mod(S')$} is either the trivial group or contains a pseudo-Anosov element. One can then construct an element $f \in H$ such that, for every connected component $S-\bigcup_{i=1}^k C_i$ such that the restriction $H|_{S'} \subseteq \Mod(S')$ contains a pseudo-Anosov element, the element \mbox{$f|_{S'} \in \Mod(S')$} is a pseudo-Anosov.

In the context of $\Out(F_{\n})$, Clay and Uyanik~\cite{clay2019atoroidal} proved Theorem~\ref{Theo intro 1} when $H$ is a subgroup of $\Out(F_{\n})$ such that $\Poly(H)=\{1\}$. Indeed, by a result of Levitt~\cite[Proposition~1.4, Lemma~1.5]{Levitt09}, if $\phi \in \Out(F_n)$ and if $\Poly(\phi) \neq\{1\}$, then there exists a nontrivial element $g \in F_{\n}$ and $k \in \NN^*$ such that $\phi^k([g])=[g]$. In this context, Clay and Uyanik proved that, if $H$ does not virtually preserve the conjugacy class of a nontrivial element of $F_{\n}$, there exists an element $\phi \in H$ which is \emph{atoroidal}: no power of $\phi$ fixes the conjugacy class of a nontrivial element of $F_{\n}$. From Clay and Uyanik's theorem, one can then ask the following question. If $H$ is a subgroup of $\Out(F_{\n})$ such that $H$ virtually fixes the conjugacy class of a nontrivial subgroup $A$ of $F_{\n}$, is it true that either $H$ virtually fixes the conjugacy class of a nontrivial element $g \in F_{\n}$ such that $g$ is not contained in a conjugate of $A$, or there exists $\phi \in H$ such that the only conjugacy classes of elements of $F_{\n}$ virtually fixed by $\phi$ are contained in a conjugate of $A$? 

Unfortunately, such a result is not true. Indeed, let $F_3=\left\langle a,b,c \right\rangle$ be a nonabelian free group of rank $3$. Let $\phi_a$ (resp. $\phi_b$) be the automorphism of $F_3$ which fixes $a$ and $b$ and which sends $c$ to $ca$ (resp. $c$ to $cb$), and let $H=\left\langle [\phi_a],[\phi_b] \right\rangle \subseteq \Out(F_3)$. Then every element $\phi \in H$ has a representative which fixes $\left\langle a,b\right\rangle$ and sends $c$ to $cg_{\phi}$ with $g_{\phi} \in \left\langle a,b \right\rangle$. Thus, $\phi$ fixes the conjugacy class of $g_{\phi}cg_{\phi}c^{-1}$. However, there always exist $\phi' \in H$, such that $\phi'$ does not preserve the conjugacy class of $g_{\phi}cg_{\phi}c^{-1}$. This example illustrates the main difficulty which appears when generalizing Clay and Uyanik's theorem: the fact that $\Poly(H) \neq \{1\}$ implies that every element of $H$ has periodic conjugacy classes which might not be fixed by the whole group. However, for the above example, we have $\Poly(H)=F_3$ and every element of $H$ satisfies Theorem~\ref{Theo intro 1}. Therefore, Theorem~\ref{Theo intro 1} is, from this viewpoint, the right generalization of Clay and Uyanik's theorem.

We now sketch the proof of Theorem~\ref{Theo intro 1}. It is inspired by the proof of \cite[Theorem~A]{clay2019atoroidal}. However, technical difficulties emerge due to the presence of elements of $F_{\n}$ with polynomial growth under iteration of elements of the considered subgroup of $\Out(F_{\n})$. The main difficulties are dealt with in the second article of the series~\cite{Guerch2021NorthSouth}. Let $H$ be a subgroup of $\Out(F_{\n})$. We first consider $H$-invariant \emph{free factor systems} $\mathcal{F}$ of $F_{\n}$, that is, $\mathcal{F}=\{[A_1],\ldots,[A_k]\}$, where, for every $i \in \{1\ldots,k\}$, $[A_i]$ is the conjugacy class of a subgroup $A_i$ of $F_{\n}$ and there exists a subgroup $B$ of $F_{\n}$ such that $F_{\n}=A_1 \ast \ldots \ast A_k \ast B$. There exists a partial order on the set of free factor systems of $F_{\tt n}$, where $\mathcal{F}_1 \leq \mathcal{F}_2$ if for every subgroup $A_1$ of $F_{\tt n}$ such that $[A_1] \in \mathcal{F}_1$, there exists a subgroup $A_2$ of $F_{\tt n}$ such that $[A_2] \in \mathcal{F}_2$ and $A_1$ is a subgroup of $A_2$. Hence we may consider a maximal $H$-invariant sequence of free factor systems $$\varnothing=\mathcal{F}_0 \leq \mathcal{F}_1 \leq\ldots \leq \mathcal{F}_k=\{[F_{\n}]\}. $$ The proof is now by induction on $i \in \{1,\ldots,k\}$: for every $i \in \{0\ldots,k\}$, we construct an element $\phi_i \in H$ such that $\Poly(\phi_i|_{\mathcal{F}_i})=\Poly(H|_{\mathcal{F}_i})$ (we define the sense of the restrictions in Section~\ref{Section currents associated automorphism}). Let $i \in \{1,\ldots,k\}$ and suppose that we have constructed $\phi_{i-1}$. There are two cases to consider. If the extension $\mathcal{F}_{i-1} \leq \mathcal{F}_i$ is \emph{nonsporadic} (see the definition in~Section~\ref{Subsection malnormal}) then the construction of $\phi_i$ from $\phi_{i-1}$ follows from the works of Handel-Mosher~\cite{HandelMosher20}, Guirardel-Horbez~\cite{Guirardelhorbez19} and Clay-Uyanik~\cite{ClayUya2018}. If the extension $\mathcal{F}_{i-1} \leq \mathcal{F}_i$ is \emph{sporadic}, the construction of $\phi_i$ relies on the action of $H$ on some natural (compact, metrizable) space that we introduced in \cite{Guerch2021currents}. This space is called the \emph{space of currents relative to $\Poly(H|_{\mathcal{F}_{i-1}})$}, denoted by $\PCurr(F_{\n},\Poly(H|_{\mathcal{F}_{i-1}}))$. It is defined as a subspace of the space of Radon measures on a natural space $\partial^2(F_{\n},\Poly(H|_{\mathcal{F}_{i-1}}))$, the double boundary of $F_{\n}$ relative to $\Poly(H|_{\mathcal{F}_{i-1}})$ (see Section~\ref{Section relative currents} for precise definitions). In~\cite{Guerch2021NorthSouth}, we proved that the element $\phi_{i-1}$ that we have constructed acts with a \emph{North-South dynamics} on $\PCurr(F_{\n},\Poly(H|_{\mathcal{F}_{i-1}}))$: there exist two proper disjoint closed subsets of $\PCurr(F_{\n},\Poly(H|_{\mathcal{F}_{i-1}}))$ such that every point of $\PCurr(F_{\n},\Poly(H|_{\mathcal{F}_{i-1}}))$ which is not contained in these subsets converges to one of the two subsets under positive or negative iteration of $\phi_{i-1}$. This North-South dynamics result allows us, applying classical ping-pong arguments similar to the one of Tits~\cite{Tits72}, to construct the element $\phi_i \in H$ such that  $\Poly(\phi_i|_{\mathcal{F}_i})=\Poly(H|_{\mathcal{F}_i})$, which concludes the proof.

The element constructed in Theorem~\ref{Theo intro 1} is in general not unique. Indeed, when the subgroup $H$ of $\Out(F_{\n})$ is such that $\Poly(H)=\{1\}$, Clay and Uyanik~\cite[Theorem~B]{clay2019atoroidal} give necessary and sufficient conditions for $H$ to contain a nonabelian free subgroup consisting in atoroidal elements.

We now outline some consequences of Theorem~\ref{Theo intro 1}. The first one is a result concerning the periodic subset of a subgroup of $\Out(F_{\n})$. Let $H$ be a subgroup of $\Out(F_{\n})$. We denote by $\mathrm{Per}(H)$ the set of conjugacy classes of $F_{\n}$ fixed by a power of every element of $H$. In the above example, we constructed a subgroup $H$ of $\Out(F_{\n})$ such that $\mathrm{Per}(H)$ contains the conjugacy class of a nonabelian subgroup of rank $2$. This is in fact the lowest possible rank where a generalization of the theorem of Clay and Uyanik using $\mathrm{Per}(H)$ instead of $\Poly(H)$ cannot work, as shown by the following result.

\begin{theo}[Corollary~\ref{Coro alternative single elements}]\label{Theo intro 2}
Let $\n \geq 3$ and let $g_1,\ldots,g_k$ be nontrivial root-free elements of $F_{\n}$. Let $H$ be subgroup of $\Out(F_{\n})$ such that, for every $i \in \{1,\ldots,k\}$, every element of $H$ has a power which fixes the conjugacy class of $g_i$. Then one of the following (mutually exclusive) statements holds.

\medskip

\noindent{$(1)$ } There exists $g_{k+1} \in F_{\n}$ such that $[g_{k+1}] \notin \{[g_1],\ldots,[g_k]\}$ and whose conjugacy class is fixed by a power of every element of $H$.

\medskip

\noindent{$(2)$ } There exists $\phi \in H$ such that $\mathrm{Per}(\phi)=\{[\left\langle g_1 \right\rangle],\ldots,[\left\langle g_k \right\rangle]\}$.
\end{theo}

As we show with Corollary~\ref{Coro pseudo Anosov subgroup Mod}, Case~$(2)$ of Theorem~\ref{Theo intro 2} naturally occurs when we are working with a subgroup of a mapping class group of a compact, connected surface $S$ whose fundamental group is identified with $F_{\n}$. Finally, we give in Proposition~\ref{Coro subgroup Fn are poly} a method, using JSJ decompositions of $F_{\n}$, allowing to compute $\Poly(H)$ for subgroups $H$ of $\Out(F_{\n})$ which act by global conjugations on some subgroups of $F_{\n}$.

\medskip

{\small{\bf Acknowledgments. } I warmly thank my advisors, Camille Horbez and Frédéric Paulin, for their precious advices and for carefully reading the different versions of this article. }

\section{Preliminaries}

\subsection{Malnormal subgroup systems of $F_{\tt n}$}\label{Subsection malnormal}

Let ${\tt n}$ be an integer greater than $1$ and let $F_{\tt n}$ be a free group of rank ${\tt n}$. A \emph{subgroup system of $F_{\tt n}$} is a finite (possibly empty) set $\mathcal{A}$ whose elements are conjugacy classes of nontrivial (that is distinct from $\{1\}$) finite rank subgroups of $F_{\tt n}$. Note that a subgroup system $\mathcal{A}$ is completely determined by the set of subgroups $A$ of $F_{\tt n}$ such that $[A] \in \mathcal{A}$. There exists a partial order on the set of subgroup systems of $F_{\tt n}$, where $\mathcal{A}_1 \leq \mathcal{A}_2$ if for every subgroup $A_1$ of $F_{\tt n}$ such that $[A_1] \in \mathcal{A}_1$, there exists a subgroup $A_2$ of $F_{\tt n}$ such that $[A_2] \in \mathcal{A}_2$ and $A_1$ is a subgroup of $A_2$. In this case we say that $\mathcal{A}_2$ is an \emph{extension} of $\mathcal{A}_1$. The \emph{stabilizer in $\Out(F_{\tt n})$ of a subgroup system} $\mathcal{A}$, denoted by $\Out(F_{\tt n},\mathcal{A})$, is the set of all elements $\phi \in \Out(F_{\tt n})$ such that $\phi(\mathcal{A})=\mathcal{A}$. If $\mathcal{A}_1$ and $\mathcal{A}_2$ are two subgroup systems, we set $\Out(F_{\n},\mathcal{A}_1,\mathcal{A}_2)=\Out(F_{\n},\mathcal{A}_1) \cap \Out(F_{\n},\mathcal{A}_2)$.

Recall that a subgroup $A$ of $F_{\tt n}$ is \emph{malnormal} if for every element $x \in F_{\tt n}-A$, we have $xAx^{-1} \cap A=\{e\}$. A subgroup system $\mathcal{A}$ is said to be \emph{malnormal} if every subgroup $A$ of $F_{\tt n}$ such that $[A] \in \mathcal{A}$ is malnormal and, for all subgroups $A_1,A_2$ of $F_{\tt n}$ such that $[A_1],[A_2] \in \mathcal{A}$, if $A_1 \cap A_2$ is nontrivial then $A_1=A_2$. An element $g \in F_{\tt n}$ is \emph{$\mathcal{A}$-peripheral} (or simply \emph{peripheral} if there is no ambiguity) if it is trivial or conjugate into one of the subgroups of $\mathcal{A}$, and \emph{$\mathcal{A}$-nonperipheral} otherwise. 

An important class of examples of malnormal subgroup systems is given by the \emph{free factor systems}. A \emph{free factor system of $F_{\tt n}$} is a (possibly empty) set $\mathcal{F}$ of conjugacy classes $\{[A_1],\ldots,[A_r]\}$ of nontrivial subgroups $A_1,\ldots,A_r$ of $F_{\tt n}$ such that there exists a subgroup $B$ of $F_{\ n}$ with $F_{\tt n}=A_1 \ast \ldots \ast A_r \ast B$. An extension $\mathcal{F}_1 \leq \mathcal{F}_2=\{[A_1],\ldots,[A_k]\}$ of free factor systems is \emph{sporadic} if there exists $\ell \in \{1,\ldots,k\}$ such that, for every $j \in \{1,\ldots,k\}-\{\ell\}$, we have $[A_j] \in \mathcal{F}_1$ and if one of the following holds:

\medskip

\noindent{$(a)$ } there exist subgroups $B_1,B_2$ of $F_{\n}$ such that $[B_1],[B_2] \in \mathcal{F}_1$ and $A_{\ell}=B_1 \ast B_2$;

\medskip

\noindent{$(b)$ } there exists a subgroup $B$ of $F_{\n}$ such that $[B] \in \mathcal{F}_1$ and $A_{\ell}$ is an HNN extension of $B$ over the trivial group;

\medskip

\noindent{$(c)$ } there exists $g \in F_{\n}$ such that $\mathcal{F}_2=\mathcal{F}_1\cup\{[g]\}$ and $A_{\ell}=\left\langle g \right\rangle$. 

\medskip 

Otherwise, the extension $\mathcal{F}_1 \leq \mathcal{F}_2$ is \emph{nonsporadic}. A free factor system $\mathcal{F}$ of $F_{\n}$ is \emph{sporadic} (resp. \emph{nonsporadic}) if the extension $\mathcal{F} \leq \{[F_{\n}]\}$ is sporadic (resp. nonsporadic). An ascending sequence of free factor systems $\mathcal{F}_1 \leq \ldots \leq \mathcal{F}_i=\{[F_{\tt n}]\}$ of $F_{\tt n}$ is called a \emph{filtration of $F_{\tt n}$}.

Given a free factor system $\mathcal{F}$ of $F_{\tt n}$, a \emph{free factor of $(F_{\tt n},\mathcal{F})$} is a subgroup $A$ of $F_{\tt n}$ such that there exists a free factor system $\mathcal{F}'$ of $F_{\tt n}$ with $[A] \in \mathcal{F}'$ and $\mathcal{F} \leq \mathcal{F}'$. When $\mathcal{F}=\varnothing$, we say that $A$ is a \emph{free factor of $F_{\tt n}$}. A free factor of $(F_{\tt n},\mathcal{F})$ is \emph{proper} if it is nontrivial, not equal to $F_{\tt n}$ and if its conjugacy class does not belong to $\mathcal{F}$. 

In general, we will work in a finite index subgroup of $\Out(F_{\n})$ defined as follows. Let $$\IA_{\n}(\ZZ/3\ZZ)=\ker(\Out(F_{\n}) \to \Aut(H_1(F_{\n},\ZZ/3\ZZ)).$$  For every $ \phi \in \IA_{\n}(\ZZ/3\ZZ)$, we have the following properties:

\medskip

\noindent{$(1)$ } any $\phi$-periodic conjugacy class of free factor of $F_{\n}$ is fixed by $\phi$ \cite[Theorem~II.3.1]{HandelMosher20};

\medskip

\noindent{$(2)$ } any $\phi$-periodic conjugacy class of elements of $F_{\n}$ is fixed by $\phi$ \cite[Theorem~II.4.1]{HandelMosher20}.

\medskip

Another class of examples of malnormal subgroup systems is the following one. Let $g \in F_{\n}$ and let $\mathfrak{B}$ be a free basis of $F_{\n}$. The length of the conjugacy class of $g$ with respect to $\mathfrak{B}$ is $$\ell_{\mathfrak{B}}([g])=\min_{h \in [g]} \ell_{\mathfrak{B}}(h),$$ where $\ell_{\mathfrak{B}}(h)$ is the word length of $h$ with respect to the basis $\mathfrak{B}$. An outer automorphism $\phi \in \Out(F_{\tt n})$ is \emph{exponentially growing} if there exists $g \in F_{\tt n}$ such that the length of the conjugacy class $[g]$ of $g$ in $F_{\tt n}$ with respect to some basis of $F_{\tt n}$ grows exponentially fast under positive iteration of $\phi$. One can show that if $g$ is exponentially growing with respect to some free basis of $F_{\n}$, then it is exponentially growing for every free basis of $F_{\n}$. If $\phi \in \Out(F_{\tt n})$ is not exponentially growing, one can show, using for instance the technology of train tracks due to Bestvina and Handel~\cite{BesHan92}, that for every $g \in F_{\n}$, the element $g$ has polynomial growth under positive iteration of $\phi$. In this case, we say that $\phi$ is \emph{polynomially growing}. A result of Levitt~\cite[Proposition~1.4~$(1)$]{Levitt09} shows that this definition is equivalent to the definition given in the introduction. For an automorphism $\alpha \in \Aut(F_{\tt n})$, we say that $\alpha$ is exponentially growing if there exists $g \in F_{\tt n}$ such that the word length of $g$ grows exponentially fast under iteration of $\phi$. Otherwise, $\alpha$ is polynomially growing. Let $\phi \in \Out(F_{\tt n})$ be exponentially growing. A subgroup $P$ of $F_{\tt n}$ is a \emph{polynomial subgroup} of $\phi$ if there exist $k \in \NN^*$ and a representative $\alpha$ of $\phi^k$ such that $\alpha(P)=P$ and $\alpha|_P$ is polynomially growing. By~\cite[Proposition~1.4]{Levitt09}, there exist finitely many conjugacy classes $[H_1],\ldots,[H_k]$ of maximal polynomial subgroups of $\phi$. Moreover, the proof of~\cite[Proposition~1.4]{Levitt09} implies that the set $\mathcal{H}=\{[H_1],\ldots,[H_k]\}$ is a malnormal subgroup system (see~\cite[Section~2.1]{Guerch2021NorthSouth}). We denote this malnormal subgroup system by $\mathcal{A}(\phi)$. Note that, if $H$ is a subgroup of $F_{\tt n}$ such that $[H] \in \mathcal{A}(\phi)$, there exists $\Phi^{-1} \in \phi^{-1}$ such that $\Phi^{-1}(H)=H$ and $\Phi^{-1}|_H$ is polynomially growing. Hence we have $\mathcal{A}(\phi) \leq \mathcal{A}(\phi^{-1})$. By symmetry, we have
\begin{equation}\label{Equation p6}
\mathcal{A}(\phi)=\mathcal{A}(\phi^{-1}).
\end{equation}
Moreover, for every element $\psi \in \Out(F_{\n})$, we have $$\mathcal{A}(\psi\phi\psi^{-1})=\psi(\mathcal{A}(\phi)).$$
In order to distinguish between the set of elements of $F_{\n}$ which have polynomial growth under positive iteration of $\phi$ and the associated malnormal subgroup system, we will denote by $\Poly(\phi)$ the former. We have $\Poly(\phi)=\Poly(\phi^{-1})$ by Equation~\eqref{Equation p6}. If $H$ is a subgroup of $\Out(F_{\n})$, we set $\Poly(H)=\bigcap_{\phi \in H} \Poly(\phi)$.

Let $\mathcal{A}$ be a malnormal subgroup system and let $\phi \in \Out(F_{\tt n},\mathcal{A})$ be a relative outer automorphism. We say that $\phi$ is \emph{atoroidal relative to $\mathcal{A}$} if, for every $k \in \NN^*$, the element $\phi^k$ does not preserve the conjugacy class of any $\mathcal{A}$-nonperipheral element. We say that $\phi$ is \emph{expanding relative to $\mathcal{A}$} if $\mathcal{A}(\phi) \leq \mathcal{A}$. Note that an expanding outer automorphism relative to $\mathcal{A}$ is in particular atoroidal relative to $\mathcal{A}$. When $\mathcal{A}=\varnothing$, then the outer automorphism $\phi$ is expanding relative to $\mathcal{A}$ if and only if for every nontrivial element $g \in F_{\tt n}$, the length of the conjugacy class $[g]$ of $g$ in $F_{\tt n}$ with respect to some basis of $F_{\tt n}$ grows exponentially fast under iteration of $\phi$. Therefore, by a result of Levitt \cite[Corollary~1.6]{Levitt09}, the outer automorphism $\phi$ is expanding relative to $\mathcal{A}=\varnothing$ if and only if $\phi$ is atoroidal relative to $\mathcal{A}=\varnothing$.

Let $\mathcal{A}=\{[A_1],\ldots,[A_r]\}$ be a malnormal subgroup system and let $\mathcal{F}$ be a free factor system. Let $i \in \{1,\ldots,r\}$. By~\cite[Theorem~3.14]{ScoWal79} for the action of $A_i$ on one of its Cayley graphs, there exist finitely many subgroups $A_i^{(1)},\ldots,A_i^{(k_i)}$ of $A_i$ such that:

\medskip

\noindent{$(1)$ } for every $j \in \{1,\ldots,k_i\}$, there exists a subgroup $B$ of $F_{\tt n}$ such that $[B] \in \mathcal{F}$ and $A_i^{(j)}=B \cap A_i$;

\medskip

\noindent{$(2)$ } for every subgroup $B$ of $F_{\tt n}$ such that $[B] \in \mathcal{F}$ and $B \cap A_i \neq \{e\}$, there exists $j \in \{1,\ldots,k_i\}$ such that $A_i^{(j)}=B \cap A_i$;

\medskip

\noindent{$(3)$ } the subgroup $A_i^{(1)} \ast \ldots \ast A_i^{(k_i)}$ is a free factor of $A_i$.

\bigskip

Thus, one can define a new subgroup system as $$\mathcal{F} \wedge \mathcal{A}=\bigcup_{i=1}^r\{[A_i^{(1)}],\ldots,[A_i^{(k_i)}]\}.$$ Since $\mathcal{A}$ is malnormal, and since, for every $i \in \{1,\ldots,r\}$, the group $A_i^{(1)} \ast \ldots \ast A_i^{(k_i)}$ is a free factor of $A_i$, it follows that the subgroup system $\mathcal{F} \wedge \mathcal{A}$ is a malnormal subgroup system of $F_{\tt n}$. We call it the \emph{meet of $\mathcal{F}$ and $\mathcal{A}$}. If $\phi \in \Out(F_{\n},\mathcal{F},\mathcal{A})$ then $\phi \in \Out(F_{\n},\mathcal{F}\wedge\mathcal{A})$.

\subsection{Relative currents}\label{Section relative currents}

In this section, we define the notion of \emph{currents of $F_{\tt n}$ relative to a malnormal subgroup system $\mathcal{A}$}. The section follows~\cite{Guerch2021currents,Guerch2021NorthSouth} (see the work of Gupta~\cite{gupta2017relative} for the particular case of free factor systems and Guirardel and Horbez \cite{Guirardelhorbez19laminations} in the context of free products of groups). It can be thought of as a functional space in which densely live the $\mathcal{A}$-nonperipheral elements of $F_{\tt n}$. 

Let $\partial_{\infty}F_{\tt n}$ be the Gromov boundary of $F_{\tt n}$. The \emph{double boundary of $F_{\tt n}$} is the Hausdorff locally compact, totally disconnected quotient topological space $$\partial^2F_{\tt n}=\left(\partial_{\infty} F_{\tt n} \times \partial_{\infty} F_{\tt n} \setminus \Delta \right)/\sim,$$ where $\sim$ is the equivalence relation generated by the flip relation $(x,y)\sim(y,x)$ and $\Delta$ is the diagonal, endowed with the diagonal action of $F_{\tt n}$. We denote by $\{x,y\}$ the equivalence class of $(x,y)$.

Let $T$ be the Cayley graph of $F_{\tt n}$ with respect to a free basis $\mathfrak{B}$. The boundary of $T$ is naturally homeomorphic to $\partial_{\infty}F_{\tt n}$ and the set $\partial^2F_{\tt n}$ is then identified with the set of unoriented bi-infinite geodesics in $T$. Let $\gamma$ be a finite geodesic path in $T$. The path $\gamma$ determines a subset in $\partial^2F_{\tt n}$ called the \emph{cylinder set of $\gamma$}, denoted by $C(\gamma)$, which consists in all unoriented bi-infinite geodesics in $T$ that contains $\gamma$. Such cylinder sets form a basis for a topology on $\partial^2 F_{\tt n}$, and in this topology, the cylinder sets are both open and closed, hence compact. The action of $F_{\tt n}$ on $\partial^2F_{\tt n}$ has a dense orbit.

Let $A$ be a nontrivial subgroup of $F_{\tt n}$ of finite rank. The induced $A$-equivariant inclusion $\partial_{\infty} A \hookrightarrow \partial_{\infty} F_{\tt n}$ induces an inclusion $\partial^2 A \hookrightarrow \partial^2 F_{\tt n}$. Let $\mathcal{A}=\{[A_1],\ldots,[A_r]\}$ be a malnormal subgroup system. Let $$\partial^2\mathcal{A}= \bigcup_{i=1}^r \bigcup_{g \in F_{\tt n}} \partial^2 \left(gA_ig^{-1}\right).$$ Let $\partial^2(F_{\tt n},\mathcal{A})=\partial^2F_{\tt n} -\partial^2\mathcal{A}$ be the \emph{double boundary of $F_{\tt n}$ relative to $\mathcal{A}$}. This subset is invariant under the action of $F_{\tt n}$ on $\partial^2F_{\tt n}$ and inherits the subspace topology of $\partial^2F_{\tt n}$.

\begin{lem}\cite[Lemmas~2.5, 2.6, 2.7]{Guerch2021currents}\label{Lem Properties of relative boundary}
Let ${\tt n} \geq 3$ and let $\mathcal{A}$ be a malnormal subgroup system of $F_{\tt n}$. The space $\partial^2(F_{\tt n},\mathcal{A})$ is an open subspace of $\partial^2 F_{\n}$, hence is locally compact, and the action of $F_{\tt n}$ on $\partial^2(F_{\tt n},\mathcal{A})$ has a dense orbit.
\end{lem}

We can now define a \emph{relative current}. Let ${\tt n} \geq 3$ and let $\mathcal{A}$ be a malnormal subgroup system of $F_{\tt n}$. A \emph{relative current on $(F_{\tt n},\mathcal{A})$} is a (possibly zero) $F_{\tt n}$-invariant nonnegative Radon measure $\mu$ on $\partial^2(F_{\tt n},\mathcal{A})$. The set $\Curr(F_{\tt n},\mathcal{A})$ of all 
relative currents on $(F_{\tt n},\mathcal{A})$ is equipped with the weak-$\ast$ topology: a sequence $(\mu_n)_{n \in \NN}$ in $\Curr(F_{\tt n},\mathcal{A})^{\NN}$ converges to a current \mbox{$\mu \in \Curr(F_{\tt n},\mathcal{A})$} if and only if for any pair of disjoint clopen subsets $S,S' \subseteq \partial^2(F_{\tt n},\mathcal{A})$, the sequence $(\mu_n(S \times S'))_{n \in \NN}$ converges to $\mu(S \times S')$.

The group $\Out(F_{\tt n},\mathcal{A})$ acts on $\Curr(F_{\tt n},\mathcal{A})$ as follows. Let $\phi \in \Out(F_{\tt n},\mathcal{A})$ and let $\Phi$ be a representative of $\phi$. The automorphism $\Phi$ acts diagonally by homeomorphisms on $\partial^2F_{\n}$. If $\Phi' \in \phi$, then the action of $\Phi'$ on $\partial^2F_{\n}$ differs from the action of $\Phi$ by a translation by an element of $F_{\n}$. Let $\mu \in \Curr(F_{\tt n},\mathcal{A})$ and let $C$ be a Borel subset of $\partial^2(F_{\tt n},\mathcal{A})$. Then, since $\phi$ preserves $\mathcal{A}$, we see that $\Phi^{-1}(C) \in \partial^2(F_{\tt n},\mathcal{A})$. Then we set $$\phi(\mu)(C)=\mu(\Phi^{-1}(C)),$$ which is well-defined since $\mu$ is $F_{\tt n}$-invariant. 

\bigskip

Every conjugacy class of nonperipheral element $g \in F_{\tt n}$ determines a relative current $\eta_{[g]}$ as follows. Suppose first that $g$ is \emph{root-free}, that is $g$ is not a proper power of any element in $F_{\tt n}$. Let $\gamma$ be a finite geodesic path in the Cayley graph $T$. Then $\eta_{[g]}(C(\gamma))$ is the number of axes in $T$ of conjugates of $g$ that contain the path $\gamma$. By~\cite[Lemma~3.2]{Guerch2021currents}, $\eta_{[g]}$ extends uniquely to a current in $\Curr(F_{\n},\mathcal{A})$ which we still denote by $\eta_{[g]}$. If $g=h^k$ with $k \geq 2$ and $h$ root-free, we set $\eta_{[g]}=k \;\eta_{[h]}$. Such currents are called \emph{rational currents}.

\bigskip

Let $\mu \in \Curr(F_{\tt n},\mathcal{A})$. The \emph{support of $\mu$}, denoted by $\Supp(\mu)$, is the support of the Borel measure $\mu$ on $\partial^2(F_{\tt n},\mathcal{A})$. We recall that $\Supp(\mu)$ is a \emph{lamination of $\partial^2(F_{\tt n},\mathcal{A})$}, that is, a closed $F_{\n}$-invariant subset of $\partial^2(F_{\tt n},\mathcal{A})$.

In the rest of the article, rather than considering the space of relative currents itself, we will consider the set of \emph{projectivized relative currents}, denoted by $$\mathbb{P}\Curr(F_{\tt n},\mathcal{A})=(\Curr(F_{\tt n},\mathcal{A})-\{0\})/\sim,$$ where $\mu \sim \nu$ if there exists $\lambda \in \RR_+^*$ such that $\mu=\lambda \nu$. The projective class of a current $\mu \in \Curr(F_{\tt n},\mathcal{A})$ will be denoted by $[\mu]$. For every $\phi \in \Out(F_{\n},\mathcal{A})$, the action $\phi\colon \mu \mapsto \phi(\mu)$ is positively linear. Therefore, the action of $\Out(F_{\n},\mathcal{A})$ on $\Curr(F_{\n},\mathcal{A})$ induces an action on $\PCurr(F_{\n},\mathcal{A})$. We have the following properties.

\begin{lem}\cite[Lemma~3.3]{Guerch2021currents}\label{Lem PCurr compact}
Let ${\tt n} \geq 3$ and let $\mathcal{A}$ be a malnormal subgroup system of $F_{\tt n}$. The space $\PCurr(F_{\tt n},\mathcal{A})$ is compact.
\end{lem}

\begin{prop}\cite[Theorem~1.2]{Guerch2021currents}\label{Prop density rational currents}
Let ${\tt n} \geq 3$ and let $\mathcal{A}$ be a malnormal subgroup system of $F_{\tt n}$. The set of projectivised rational currents about nonperipheral elements of $F_{\tt n}$ is dense in $\PCurr(F_{\tt n},\mathcal{A})$.
\end{prop}

\subsection{Currents associated with an almost atoroidal outer automorphism of $F_{\tt n}$}\label{Section currents associated automorphism}

Let ${\tt n} \geq 3$ and let $\mathcal{F}=\{[A_1],\ldots,[A_k]\}$ be a free factor system of $F_{\tt n}$. If $\phi \in \IA_{\n}(\ZZ/3\ZZ)$ preserves $\mathcal{F}$, we denote by $$\phi|_{\mathcal{F}}=([\Phi_1|_{A_1}],\ldots,[\Phi_k|_{A_k}]) \in \prod_{i=1}^k \Out(A_i)$$ where, for every $i \in \{1,\ldots,k\}$, the element $\Phi_i$ is a representative of $\phi$ such that $\Phi_i(A_i)=A_i$. Note that the outer class of $\Phi_i|_{A_i}$ in $\Out(A_i)$ does not depend on the choice of $\Phi_i$ since $A_i$ is a malnormal subgroup of $F_{\n}$. Note that, for every $i \in \{1,\ldots,k\}$, the element $[\Phi_i|_{A_i}]$ is expanding relative to $\mathcal{F} \wedge \{[A_i]\}=\{[A_i]\}$. Hence we will say that \emph{$\phi|_{\mathcal{F}}$ is expanding relative to $\mathcal{F}$}. Let $$\Poly(\phi|_{\mathcal{F}})=\bigcup_{i=1}^k\bigcup_{g \in F_{\n}} g\Poly([\Phi_i|_{A_i}])g^{-1} \subseteq F_{\n}.$$ If $H$ is a subgroup of $\IA_{\n}(\ZZ/3\ZZ)$ which preserves $\mathcal{F}$, we set $\Poly(H|_\mathcal{F})=\bigcap_{\phi \in H} \Poly(\phi|_{\mathcal{F}})$.

We now define a class of outer automorphisms of $F_{\tt n}$ which we will study in the rest of the article.

\begin{defi}\label{Defi almost atoroidal outer automorphism}
Let ${\tt n} \geq 3$ and let $\mathcal{F}$ be a free factor system of $F_{\tt n}$. Let $\phi \in \IA_{\n}(\ZZ/3\ZZ)$. The outer automorphism $\phi$ is \emph{almost atoroidal relative to} $\mathcal{F}$ if $\Poly(\phi) \neq \{[F_{\n}]\}$ and if one of the following holds:

\medskip

\noindent{$(1)$ } $\phi$ is an atoroidal outer automorphism relative to $\mathcal{F}$; 

\medskip

\noindent{$(2)$ } the extension $\mathcal{F} \leq \{[F_{\n}]\}$ is sporadic.
\end{defi}

Let $\phi \in \IA_{\n}(\ZZ/3\ZZ)$ be an almost atoroidal outer automorphism relative to $\mathcal{F}$. In this section, we recall from~\cite{Guerch2021NorthSouth} the definition and some properties of some subsets of $\PCurr(F_{\tt n},\mathcal{F} \wedge \mathcal{A}(\phi))$ associated with $\phi$. Let $K_{PG}(\phi)$ be the subspace of all currents in  $\PCurr(F_{\tt n},\mathcal{F} \wedge \mathcal{A}(\phi))$ whose support is contained in $\partial^2\mathcal{A}(\phi) \cap \partial^2(F_{\tt n},\mathcal{F} \wedge \mathcal{A}(\phi))$. The subspace $K_{PG}(\phi)$ is called the \emph{space of polynomially growing currents associated with $\phi$}.

\begin{prop}\cite[Proposition~4.4, Proposition~4.12, Proposition~5.23]{Guerch2021NorthSouth}\label{Prop Existence and properties of Delta}
Let ${\tt n} \geq 3$ and let $\mathcal{F}$ be a free factor system of $F_{\tt n}$. Let $\phi \in \IA_{\n}(\ZZ/3\ZZ)$ be an almost atoroidal outer automorphism relative to $\mathcal{F}$. There exist two unique proper compact $\phi$-invariant subsets $\Delta_{\pm}(\phi)$ of $\PCurr(F_{\tt n},\mathcal{F \wedge }\mathcal{A}(\phi))$ such that the following holds. For every $[\mu] \in \Delta_+(\phi) \cup \Delta_-(\phi)$, the support of $\mu$ is contained in $\partial^2\mathcal{F}$. Let $U_+$ be a neighborhood of $\Delta_+(\phi)$, let $U_-$ be a neighborhood of $\Delta_-(\phi)$, let $V$ be a neighborhood of $K_{PG}(\phi)$. There exists $N \in \NN^*$ such that for every $n \geq 1$ and every ($\mathcal{F} \wedge \mathcal{A}(\phi)$)-nonperipheral $w \in F_{\tt n}$ such that $\eta_{[w]} \notin V$, one of the following holds $$\phi^{Nn}(\eta_{[w]}) \in U_+ \text{~~~~ or ~~~~} \phi^{-Nn}(\eta_{[w]}) \in U_-.$$
\end{prop}

The subsets $\Delta_+(\phi)$ and $\Delta_-(\phi)$ are called the \emph{simplices of attraction and repulsion of $\phi$}. Let $\mathcal{F} \leq \mathcal{F}_1=\{[A_1],\ldots,[A_k]\}$ be two free factor systems of $F_{\n}$. Let $\phi \in \IA_{\n}(\ZZ/3\ZZ) \cap \Out(F_n,\mathcal{F},\mathcal{F}_1)$. We say that  $\phi|_{\mathcal{F}_1}$ is \emph{almost atoroidal relative to $\mathcal{F}$} if, for every $i \in \{1,\ldots,k\}$, the outer automorphism $[\Phi_i|_{A_i}]$ is almost atoroidal relative to $\mathcal{F} \wedge \{[A_i]\}$. Let $i \in \{1,\ldots,k\}$. If $\phi|_{\mathcal{F}_1}$ is almost atoroidal relative to $\mathcal{F}$, we denote by $\Delta_{\pm}([A_i],\phi) \subseteq \PCurr(A_i,\mathcal{F} \wedge \{[A_i]\} \wedge \mathcal{A}([\Phi_i|_{A_i}]))$ the convexes of attraction and repulsion of $[\Phi_i|_{A_i}]$. If $\psi \in \IA_{\n}(\ZZ/3\ZZ)$ preserves the conjugacy class of $A_i$ and $\mathcal{F} \wedge \{[A_i]\} \wedge \mathcal{A}([\Phi_i|_{A_i}])$, then $\Delta_{\pm}([A_i],\psi\phi\psi^{-1})=\psi(\Delta_{\pm}([A_i],\phi))$.

\bigskip

We will also need the following result which gives the existence and properties of an approximation of the length function of the conjugacy class of an element if $F_{\n}$ in the context of the space of currents.

\begin{prop}\cite[Lemma~3.26, Lemma~3.27~$(3)$]{Guerch2021NorthSouth}\label{Prop moving exponential1}
Let ${\tt n} \geq 3$ and let $\mathcal{F}$ be a free factor system of $F_{\tt n}$. Let $\phi \in \Out(F_{\tt n},\mathcal{F})$ be an almost atoroidal outer automorphism relative to $\mathcal{F}$ of type~$(2)$. There exists a continuous, positively linear function $\lVert. \rVert_{\mathcal{F}} \colon \Curr(F_n,\mathcal{F} \wedge \mathcal{A}(\phi)) \to \RR_+$ such that the following holds.

\medskip

\noindent{$(1)$ } There exist a basis $\mathfrak{B}$ of $F_{\n}$ and a constant $C>1$ such that, for every $\mathcal{F} \wedge \mathcal{A}(\phi)$-nonperipheral element $g \in F_{\n}$, we have $\lVert \eta_{[g]} \rVert_{\mathcal{F}} \in \NN^*$ and $$\ell_{\mathfrak{B}}([g]) \geq C\;\lVert \eta_{[g]} \rVert_{\mathcal{F}}.$$

\end{prop}

\noindent{$(2)$ } For every $\eta \in \Curr(F_{\n},\mathcal{F}\wedge\mathcal{A}(\phi))$, if $\lVert \eta \rVert_{\mathcal{F}}=0$, then $\eta=0$.

\medskip

Let $$\widehat{\Delta}_{\pm}(\phi)=\left\{[t\mu +(1-t)\nu]\;\left|\; t \in [0,1], [\mu] \in \Delta_{\pm}(\phi), [\nu] \in K_{PG}(\phi), \lVert \mu \rVert_{\mathcal{F}}=\lVert \nu \rVert_{\mathcal{F}}=1 \right.\right\}$$ be the \emph{convexes of attraction and repulsion of $\phi$}. We have the following result. 

\begin{theo}\cite[Theorem~6.4]{Guerch2021NorthSouth}\label{Theo North-South dynamics almost atoroidal}
Let ${\tt n} \geq 3$ and let $\mathcal{F}$ be a free factor system of $F_{\tt n}$. Let $\phi \in \IA_{\n}(\ZZ/3\ZZ)\cap \Out(F_{\tt n},\mathcal{F})$ be an almost atoroidal outer automorphism relative to $\mathcal{F}$ of type~$(2)$. Let $\widehat{\Delta}_{\pm}(\phi)$ be the convexes of attraction and repulsion of $\phi$ and $\Delta_{\pm}(\phi)$ be the simplices of attraction and repulsion of $\phi$. Let $U_{\pm}$ be open neighborhoods of $\Delta_{\pm}(\phi)$ in $\PCurr(F_{\tt n},\mathcal{F} \wedge \mathcal{A}(\phi))$ and $\widehat{V}_{\pm}$ be open neighborhoods of $\widehat{\Delta}_{\pm}(\phi)$ in $\PCurr(F_{\tt n},\mathcal{F} \wedge \mathcal{A}(\phi))$. There exists $M \in \NN^*$ such that for every $n \geq M$, we have 
$$
\phi^{\pm n}(\PCurr(F_{\tt n},\mathcal{F} \wedge \mathcal{A}(\phi))-\widehat{V}_{\mp}) \subseteq U_{\pm}.
$$
\end{theo}

\begin{prop}\cite[Corollary~6.5]{Guerch2021NorthSouth}\label{Prop moving exponential}
Let ${\tt n} \geq 3$ and let $\mathcal{F}$ be a free factor system of $F_{\tt n}$. Let $\phi \in \Out(F_{\tt n},\mathcal{F})$ be an almost atoroidal outer automorphism relative to $\mathcal{F}$ of type~$(2)$. There exists a continuous, positively linear function $\lVert. \rVert_{\mathcal{F}} \colon \Curr(F_n,\mathcal{F} \wedge \mathcal{A}(\phi)) \to \RR_+$ such that the following holds.

For every open neighborhood $\widehat{V}_{-} \subseteq \PCurr(F_{\tt n},\mathcal{F} \wedge \mathcal{A}(\phi))$ of $\widehat{\Delta}_-(\phi)$, there exists $M \in \NN^*$ and a constant $L_0>0$ such that, for every current $[\mu] \in \PCurr(F_{\tt n},\mathcal{F} \wedge \mathcal{A}(\phi))-\widehat{V}_-$, and every $m \geq M$, we have $$\lVert\phi^m(\mu) \rVert_{\mathcal{F}} \geq 3^{m-M}L_0 \lVert\mu \rVert_{\mathcal{F}}.$$ 
\end{prop}

\section{Nonsporadic extensions and fully irreducible outer automorphisms}

Let $\n \geq 3$ and let $\mathcal{F}$ and $\mathcal{F}_1=\{[A_1],\ldots,[A_k]\}$ be two free factor systems of $F_{\n}$ with $\mathcal{F} \leq \mathcal{F}_1$ such that the extension $\mathcal{F} \leq \mathcal{F}_1$ is nonsporadic. Let $H$ be a subgroup of $\IA_{\n}(\ZZ/3\ZZ)$ which preserves $\mathcal{F}$ and $\mathcal{F}_1$. We suppose that $H$ is \emph{irreducible with respect to $\mathcal{F} \leq \mathcal{F}_1$}, that is, there does not exist a proper, nontrivial free factor system $\mathcal{F}'$ of $F_{\n}$ preserved by $H$ with $\mathcal{F} < \mathcal{F}' < \mathcal{F}_1$. Suppose that there exists $\phi \in H$ such that $\Poly(\phi|_{\mathcal{F}})=\Poly(H|_{\mathcal{F}})$. In this section, we show that there exists $\widehat{\phi} \in H$ such that $\Poly(\widehat{\phi}|_{\mathcal{F}_1})=\Poly(H|_{\mathcal{F}_1})$. The key point is to construct \emph{fully irreducible outer automorphisms relative to $\mathcal{F}$} in $H$ in the following sense. Let $\phi \in \Out(F_{\n},\mathcal{F})$. We say that $\phi$ is \emph{fully irreducible relative to $\mathcal{F}$} if no power of $\phi$ preserves a proper free factor system $\mathcal{F}'$ of $F_{\tt n}$ such that $\mathcal{F} < \mathcal{F}'$. If $\phi \in \Out(F_{\n},\mathcal{F},\mathcal{F}_1)$, we say that $\phi|_{\mathcal{F}_1}$ is \emph{fully irreducible relative to $\mathcal{F}$} (resp. \emph{atoroidal relative to $\mathcal{F}$}) if, for every $i \in \{1,\ldots,k\}$, the outer automorphism $[\Phi_i|_{A_i}]$ is fully irreducible relative to $\mathcal{F} \wedge \{[A_i]\}$ (resp. atoroidal relative to $\mathcal{F} \wedge \{[A_i]\}$). If $H$ is a subgroup of $\Out(F_n,\mathcal{F},\mathcal{F}_1)$, we say that $H|_{\mathcal{F}_1}$ is \emph{atoroidal relative to $\mathcal{F}$} if there  does not exist a conjugacy class of $F_{\n}$ which is $H$-invariant, $\mathcal{F}$-nonperiperal and $\mathcal{F}_1$-peripheral. First, we recall some properties of fully irreducible outer automorphisms.

\begin{prop}\label{Prop properties fully irreducible}
Let ${\tt n } \geq 3$ and let $\mathcal{F}$ be a nonsporadic free factor system of $F_{\tt n}$. Let $H$ be a subgroup of $\IA_{\n}(\ZZ/3\ZZ)$ which preserves $\mathcal{F}$ and such that $H$ is irreducible with respect to the extension $\mathcal{F} \leq \{[F_{\n}]\}$. Let $\phi \in H$ be a fully irreducible outer automorphism relative to $\mathcal{F}$. 

\medskip

\noindent{$(1)$ }\cite[Corollary~3.14]{Guerch2021NorthSouth} There exists at most one (up to taking inverse) conjugacy class $[g]$ of root-free $\mathcal{F}$-nonperipheral element of $F_{\tt n}$ which has polynomial growth under iteration of $\phi$. Moreover, the conjugacy class $[g]$ is fixed by $\phi$.

\medskip 

\noindent{$(2)$ }\cite[Theorem~7.4]{Guirardelhorbez19} One of the following holds:

\medskip

\noindent{$(a)$ } there exists $\psi \in H$ such that $\psi$ is a fully irreducible, atoroidal outer automorphism relative to $\mathcal{F}$; 

\medskip

\noindent{$(b)$ } if $\phi$ fixes the conjugacy class of a root-free $\mathcal{F}$-nonperipheral element $g$ of $F_{\n}$, then $[g]$ is fixed by $H$.
\end{prop}

Hence Proposition~\ref{Prop properties fully irreducible} suggests that, if $H$ is a subgroup of $F_{\n}$ which satisfies the hypotheses in Proposition~\ref{Prop properties fully irreducible}, one step in order to prove Theorem~\ref{Theo intro 1} is to construct relative fully irreducible (atoroidal) outer automorphisms in $H$. This is contained in the following theorem.

\begin{theo}\label{Theo nonsporadic extension}
Let $\n \geq 3$ and let $H$ be a subgroup of $\IA_{\n}(\ZZ/3\ZZ)$. Let $$\varnothing=\mathcal{F}_0 < \mathcal{F}_1 <\ldots < \mathcal{F}_k=\{[F_{\n}]\} $$ be a maximal $H$-invariant sequence of free factor systems. There exists $\phi \in H$ such that for every $i \in \{1,\ldots,k\}$ such that the extension $\mathcal{F}_{i-1} \leq \mathcal{F}_i$ is nonsporadic, the element $\phi|_{\mathcal{F}_i}$ is fully irreducible relative to $\mathcal{F}_{i-1}$. Moreover, if $H|_{\mathcal{F}_i}$ is atoroidal relative to $\mathcal{F}_{i-1}$, one can choose $\phi$ so that $\phi|_{\mathcal{F}_i}$ is atoroidal relative to $\mathcal{F}_{i-1}$.
\end{theo}

\dem The proof follows~\cite[Theorem~6.6]{ClayUya2018} (see also~\cite[Corollary~3.4]{clay2019atoroidal}). Let $$S=\{j\;|\; \text{the extension } \mathcal{F}_{j-1} \leq \mathcal{F}_j \text{ is nonsporadic}\}$$ and let $j \in S$. 

\medskip

\noindent{\bf Claim. } There exists a unique conjugacy class $[B_j]$ of a subgroup $B_j$ in $F_{\n}$ such that $[B_j] \in \mathcal{F}_j$ and $[B_j] \notin \mathcal{F}_{j-1}$. 

\medskip

\dem There exists at least one such conjugacy class since $\mathcal{F}_{j-1}<\mathcal{F}_j$. Suppose towards a contradiction that there exist two distinct subgroups $B_+$ and $B_-$ of $F_{\n}$ such that $[B_+] \neq [B_-]$, $[B_+],[B_-] \in \mathcal{F}_j$ and $[B_+],[B_-] \notin \mathcal{F}_{j-1}$. Then $$\mathcal{F}'([B_-])=(\mathcal{F}_j-\{[B_+]\}) \cup (\mathcal{F}_{j-1} \wedge \{[B_+]\})$$ is $H$-invariant and $\mathcal{F}_{j-1} < \mathcal{F}'([B_-])< \mathcal{F}_j$, which contradicts the maximality hypothesis of the sequence of free factor systems. The claim follows.
\hfill\qedsymbol

\medskip

Let $B_j$ be a subgroup of $F_{\n}$ given by the claim. Let $A_{j,1},\ldots,A_{j,s}$ be the subgroups of $B_j$ with pairwise disjoint conjugacy classes such that $\mathcal{A}_{j-1}=\{[A_{j,1}],\ldots,[A_{j,s}]\} \subseteq \mathcal{F}_{j-1}$ and $s$ is maximal for this property. By~\cite[Theorem~D]{HandelMosher20}, for every $j \in S$, there exists $\phi \in H$ such that $[\Phi_j|_{A_j}] \in \Out(B_j,\mathcal{A}_{j-1})$ is fully irreducible relative to $\mathcal{A}_{j-1}$. By Proposition~\ref{Prop properties fully irreducible}~$(2)$, for every $j \in S$ such that $H|_{\mathcal{F}_j}$ is atoroidal relative to $\mathcal{F}_{j-1}$, there exists $\phi \in H$ such that $[\Phi_j|_{A_j}] \in \Out(B_j,\mathcal{A}_{j-1})$ is fully irreducible and atoroidal relative to $\mathcal{A}_{j-1}$. Note that, for every $j \in S$, the free factor system $\mathcal{A}_{j-1}$ is a nonsporadic free factor system of $B_j$ by the claim and since the extension $\mathcal{F}_{j-1} \leq \mathcal{F}_j$ is nonsporadic. Let $S_1$ be the subset of $S$ consisting in every $j \in S$ such that $H|_{\mathcal{F}_j}$ is atoroidal relative to $\mathcal{F}_{j-1}$, and let $S_2=S-S_1$. By~\cite[Theorem~4.1,4.2]{Guirardelhorbez19} (see also~\cite{Mann2014,Mann2014Thesis,Horbez17cyclic,gupta18}), for every $j \in S_1$ (resp. $j \in S_2$) there exists a Gromov-hyperbolic space $X_j$ (the \emph{$\mathcal{Z}$-factor complex of $B_j$ relative to $\mathcal{A}_{j-1}$} when $j \in S_1$ and the \emph{free factor complex of $B_j$ relative to $\mathcal{A}_{j-1}$} otherwise) on which $\Out(B_j,\mathcal{A}_{j-1})$ acts by isometries and such that $\phi_0 \in \Out(B_j,\mathcal{A}_{j-1})$ is a loxodromic element if and only if $\phi_0$ is fully irreducible atoroidal relative to $\mathcal{A}_{j-1}$ (resp. fully irreducible relative to $\mathcal{A}_{j-1}$). The conclusion then follows from~\cite[Theorem~5.1]{ClayUya2018}.
\hfill\qedsymbol

\section{Sporadic extensions and polynomial growth}\label{Section sporadic extension}

Let $\n \geq 3$ and let $\mathcal{F}$ and $\mathcal{F}_1=\{[A_1],\ldots,[A_k]\}$ be two free factor systems of $F_{\n}$ with $\mathcal{F} \leq \mathcal{F}_1$. Suppose that the extension $\mathcal{F} \leq \mathcal{F}_1$ is sporadic. Let $H$ be a subgroup of $\IA_{\n}(\ZZ/3\ZZ) \cap \Out(F_{\n},\mathcal{F},\mathcal{F}_1)$. 

In order to prove Theorem~\ref{Theo intro 1}, we need to show that if $\Poly(\phi|_{\mathcal{F}})=\mathrm{Poly}(H|_{\mathcal{F}})$, there exists $\psi \in H$ such that $\Poly(\psi|_{\mathcal{F}_1})=\mathrm{Poly}(H|_{\mathcal{F}_1})$. Let $\phi \in H$ be such that $\Poly(\phi|_{\mathcal{F}})=\mathrm{Poly}(H|_{\mathcal{F}})$. Note that, for every element $g$ of $\Poly(\phi|_{\mathcal{F}})$, there exists a subgroup $A$ of $F_{\n}$ such that $[A] \in \mathcal{F} \wedge \mathcal{A}(\phi)$ and $g \in A$. Conversely, for every subgroup $A$ of $F_{\n}$ such that $[A] \in \mathcal{F} \wedge \mathcal{A}(\phi)$ and every element $g \in A$, we have $g \in \Poly(\phi|_{\mathcal{F}})$. Thus $\mathcal{F}\wedge \mathcal{A}(\phi)$ is the natural malnormal subgroup system associated with $\Poly(\phi|_{\mathcal{F}})=\Poly(H|_{\mathcal{F}})$. Thus, we see that $H$ preserves $\mathcal{F} \wedge \mathcal{A}(\phi)$ and hence $H$ acts by homeomorphisms on $\PCurr(F_n,\mathcal{F} \wedge \mathcal{A}(\phi))$. 

\begin{lem}\label{Lem Kpg moved by elements}
Let $\n \geq 3$, let $\mathcal{F}$ be a sporadic free factor system of $F_{\n}$ and let $H$ be a subgroup of $\IA_{\n}(\ZZ/3\ZZ) \cap \Out(F_{\n},\mathcal{F})$ which is irreducible with respect to $\mathcal{F} \leq \{[F_{\n}]\}$. Suppose that there exists $\phi \in H$ such that $\Poly(\phi|_{\mathcal{F}})=\mathrm{Poly}(H|_{\mathcal{F}})$. If $\Poly(\phi) \neq \Poly(H)$, there exists an infinite subset $X \subseteq H$ such that for all distinct $\psi_1,\psi_2 \in X$, we have $\psi_1(K_{PG}(\phi)) \cap \psi_2(K_{PG}(\phi))=\varnothing$.
\end{lem}

\dem Let $\mathcal{F}\wedge \mathcal{A}(\phi)=\{[A_1],\ldots,[A_r]\}$. Suppose towards a contradiction that $\mathcal{A}(\phi)=\mathcal{F} \wedge \mathcal{A}(\phi)$. Then $$\Poly(\phi)=\Poly(\phi|_{\mathcal{F}})=\Poly(H|_{\mathcal{F}})=\Poly(H).$$ This contradicts the fact that $\Poly(\phi) \neq \Poly(H)$. Thus, we have $\mathcal{A}(\phi) \neq \mathcal{F} \wedge \mathcal{A}(\phi)$. By~\cite[Lemma~5.18~$(7)$]{Guerch2021NorthSouth}, one of the following holds.

\medskip

\noindent{$(i)$ } There exist distinct $i,j \in \{1,\ldots,r\}$ such that, up to replacing $A_i$ by a conjugate, we have $\mathcal{A}(\phi)=(\mathcal{F}\wedge \mathcal{A}(\phi)-\{[A_i],[A_j]\}) \cup \{[A_i \ast A_j]\}$.

\medskip

\noindent{$(ii)$ } There exists $i \in \{1,\ldots,r\}$ and an element $g \in F_{\n}$ such that $\mathcal{A}(\phi)=(\mathcal{F}\wedge \mathcal{A}(\phi)-\{[A_i]\}) \cup \{[A_i \ast \left\langle g \right\rangle] \}$.

\medskip

\noindent{$(iii)$ } There exists $g \in F_{\n}$ such that $\mathcal{A}(\phi)=\mathcal{F}\wedge\mathcal{A}(\phi) \cup \{[\left\langle g \right\rangle]\}$.

\medskip

\noindent{\bf Case~1 } Suppose that there exist distinct $i,j \in \{1,\ldots,r\}$ such that $$\mathcal{A}(\phi)=(\mathcal{F}\wedge \mathcal{A}(\phi)-\{[A_i],[A_j]\}) \cup \{[A_i \ast A_j]\}.$$ 

\medskip

Since $\Poly(\phi|_{\mathcal{F}})=\mathrm{Poly}(H|_{\mathcal{F}})$ and $\Poly(\phi) \neq \Poly(H)$, there exists $\psi \in H$ such that, for every $n \in \NN^*$, the element $\psi^n$ does not preserve $[A_i \ast A_j]$ while preserving $[A_i]$ and $[A_j]$. Hence there exist a representative $\Psi$ of $\psi$ such that, for every $n \in \NN^*$, there exists $g_n \in  F_n-A_i \ast A_j$ such that  $\Psi^n(A_i)=A_i$ and $\Psi^n(A_j)=g_nA_jg_n^{-1}$. 

\medskip

\noindent{\bf Claim~1. } For every $n \in \NN^*$, every $g \in F_{\n}$ and every $h \in F_n$, if $$h \in (g\Psi^n(A_i \ast A_j)g^{-1}) \cap (A_i \ast A_j),$$ then $h$ is either contained in a conjugate of $A_i$ or a conjugate of $A_j$. 

\medskip

\dem Let $n \in \NN^*$ and let $h \in (g\Psi^n(A_i \ast A_j)g^{-1}) \cap (A_i \ast A_j)$. Suppose towards a contradiction that $h$ is not contained in a conjugate of $A_i$ or a conjugate of $A_j$. By~\cite[Lemma~1.2]{Levitt09}, there exists a nontrivial $\RR$-tree $T$ equipped with a minimal, isometric action of $F_{\n}$ with trivial edge stabilizers and such that every polynomial subgroup of $\phi$ fixes a point in $T$. 

The groups $A_i \ast A_j$, $gg_n A_i \ast A_jg_n^{-1}g^{-1}$ and $gA_i \ast A_jg^{-1}$ fix points in $T$. Note that, if we have $gg_n \in A_i \ast A_j$, then, since $g_n \notin A_i \ast A_j$, we have $g \notin A_i \ast A_j$. By malnormality of $A_i \ast A_j$, we have $(A_i \ast A_j) \cap (gA_i \ast A_jg^{-1})=\{1\}$. Thus, for every $g \in F_{\n}$, one of the following holds: $(A_i \ast A_j) \cap (gA_i \ast A_jg^{-1})=\{1\}$ or $(A_i \ast A_j) \cap (gg_nA_i \ast A_jg_n^{-1}g^{-1})=\{1\}$. If $A_i \ast A_j$, $gA_i \ast A_gg^{-1}$ and $gg_n (A_i \ast A_j)g_n^{-1}g^{-1}$ fix the same point $x$, then, by induction on the rank of $\Stab(x)$ (which is less than $\n$ by~\cite{GabLev95}), one can construct a nontrivial $\RR$-tree $T'$ equipped with a minimal, isometric action of a subgroup $B'$ of $F_{\n}$ containing $A_i \ast A_j$, $gA_i \ast A_jg^{-1}$ and $gg_n A_i \ast A_jg_n^{-1}g^{-1}$ with trivial arc stabilizers, such that $A_i \ast A_j$ fixes a point $x_1$, $gA_i \ast A_jg^{-1}$ fixes a point $x_2$, $gg_n A_i \ast A_jg_n^{-1}g^{-1}$ fixes a point $x_3$ and one of the following holds: $x_1 \neq x_2$ or $x_1 \neq x_3$. 

Suppose first that $x_2=x_3$. Then $x_1 \neq x_2$. Since $g\Psi^n(A_i \ast A_j)g^{-1}=gA_i \ast (g_nA_jg_n^{-1})g^{-1}$, $h$ fixes both $x_1$ and $x_2$. This contradicts the fact that $T'$ has trivial arc stabilizers.

Suppose now that $x_2 \neq x_3$. Since $g\Psi^n(A_i \ast A_j)g^{-1}=gA_i \ast (g_nA_jg_n^{-1})g^{-1}$, and since $h$ is not contained in a conjugate of $A_i$ or a conjugate of $A_j$, the element $h$ can be written as a product of elements $a_1b_1\ldots a_kb_k$ where, for every $i \in \{1,\ldots,k\}$, $a_i$ fixes $x_2$ and $b_i$ fixes $x_3$.

We claim that $h$ is loxodromic in $T'$. Let 
$G=\left\langle \Stab(x_2),\Stab(x_3)\right\rangle$. The minimal tree $T_G'$ in $T'$ of $G$ is simplicial with trivial edge stabilizers and the quotient $T_G'/G$ has exactly one edge. Hence if $g' \in G$ stabilizes a point in $T'$ it is either contained in a conjugate of $\Stab(x_2)$ or a conjugate of $\Stab(x_3)$. We may suppose that $h$ is a cyclically reduced element when written in the generating set $\{\Stab(x_2),\Stab(x_3)\}$. In particular, $h$ either fixes $x_2$ or $x_3$. Since $h$ is not contained in a conjugate of $A_i$ or a conjugate of $A_j$, we have $k \geq 2$. Hence $h$ cannot fix $x_2$ or $x_3$ and $h$ is a loxodromic element. Therefore $h$ does not fix $x_1$ and $h \notin A_i \ast A_j$, a contradiction. 
\hfill\qedsymbol

\medskip

Claim~$1$ implies that, for every distinct $m,n \in \NN^*$ and every $x \in F_{\n}$, the intersection $\Psi^n(A_i \ast A_j) \cap (x\Psi^m(A_i \ast A_j)x^{-1})$ is either contained in a conjugate of $A_i$ or a conjugate of $A_j$. By for instance \cite[Fact~I.1.2]{HandelMosher20}, for every distinct $m,n \in \NN^*$ and every $x \in F_{\n}$, we have 
$$\begin{array}{cccl}
\partial^2 \left(\Psi^n(A_i \ast A_j)\right) \cap \partial^2 \left(x\Psi^m(A_i \ast A_j)x^{-1}\right)&=&\partial^2\left( \Psi^n(A_i \ast A_j) \cap x\Psi^m(A_i \ast A_j)x^{-1}\right) \\ 
{} & \subseteq & \overline{\bigcup_{y \in F_{\n}} \big(\partial^2 \left(yA_iy^{-1}\right) \cup\partial^2 \left(yA_jy^{-1}\right)\big)}.
\end{array}$$ 
By definition of $K_{PG}(\phi)$, we have $[\mu] \in K_{PG}(\phi)$ if and only if $$\Supp(\mu) \subseteq \partial^2\mathcal{A}(\phi) \cap \partial^2(F_{\n},\mathcal{F}\wedge \mathcal{A}(\phi)) =\partial^2 \{[A_i \ast A_j]\}\cap \partial^2(F_{\n},\mathcal{F}\wedge \mathcal{A}(\phi)).$$ Moreover, if $n \in \NN^*$ and if $[\mu] \in \psi^n(K_{PG}(\phi))$, then $$\Supp(\mu) \subseteq \partial^2\psi^n(\mathcal{A}(\phi)) \cap \partial^2(F_{\n},\mathcal{F}\wedge \mathcal{A}(\phi)) =\partial^2 \{[A_i \ast g_nA_jg_n^{-1}]\}\cap \partial^2(F_{\n},\mathcal{F}\wedge \mathcal{A}(\phi)).$$ Let $n,m \in \NN^*$ be distinct. Suppose towards a contradiction that $$\psi^n(K_{PG}(\phi)) \cap \psi^m(K_{PG}(\phi)) \neq \varnothing$$ and let $[\mu] \in \psi^n(K_{PG}(\phi)) \cap \psi^m(K_{PG}(\phi))$. By $F_{\n}$-invariance of $\mu$, there exists $x \in F_{\n}$ such that $\mu$ gives positive measure to 
\begin{align*}
\partial^2 (A_i \ast g_nA_jg_n^{-1}) \cap \partial^2 \left(x(A_i \ast g_mA_jg_m^{-1})x^{-1}\right) \cap\partial^2(F_{\n},\mathcal{F}\wedge \mathcal{A}(\phi)) \\
\subseteq \overline{\Big(\bigcup_{y \in F_{\n}} \partial^2 (yA_iy^{-1}) \cup \partial^2 (yA_jy^{-1})\Big)} \cap \partial^2(F_{\n},\mathcal{F}\wedge \mathcal{A}(\phi))
\end{align*} and the last intersection is empty by the definition of the relative boundary, a contradiction.

\medskip

\noindent{\bf Case~2 } Suppose that either there exists $i \in \{1,\ldots,r\}$ and an element $g \in F_{\n}$ such that $\mathcal{A}(\phi)=(\mathcal{F}\wedge \mathcal{A}(\phi)-\{[A_i]\}) \cup \{[A_i \ast \left\langle g \right\rangle] \}$ or there exists $g \in F_{\n}$ such that $\mathcal{A}(\phi)=\mathcal{F}\wedge\mathcal{A}(\phi) \cup \{[\left\langle g \right\rangle]\}$. 

\medskip

In order to treat both cases simultaneously, in the case that there exists $g \in F_{\n}$ such that $\mathcal{A}(\phi)=\mathcal{F}\wedge\mathcal{A}(\phi) \cup \{[\left\langle g \right\rangle]\}$, we fix $A_i=\{e\}$. Case~$(2)$ only occurs when the extension $\mathcal{F} \leq \{[F_{\n}]\}$ is an HNN extension over the trivial group. In particular, we have $\mathcal{F}=\{[A]\}$ for some subgroup $A$ of $F_{\n}$ and, up to changing the representative of $[A]$, we have $F_{\n}=A \ast \left\langle g \right\rangle$ and $A_i \subseteq A$. In particular, since $H$ preserves the extension $\mathcal{F} \leq \{[F_{\n}]\}$, for every $\psi \in H$, there exist a unique representative $\Psi_0$ of $\psi$ and $g_{\psi} \in A$ such that $\Psi_0(A)=A$ and $\Psi_0(g)=gg_{\psi}$. Since $H$ is irreducible with respect to $\mathcal{F} \leq \{[F_{\n}]\}$, the subgroup $H$ does not preserve the free factor system $\mathcal{F} \cup \{[g]\}$. Thus, there exists $\psi' \in H$ such that $g_{\psi'} \neq 1$. We claim that there exist $\psi \in H$ with $g_{\psi} \notin A_i$, a representative $\Psi$ of $\psi$ and $h_1 \in A$ such that $\Psi(A_i)=A_i$ and $\Psi(g)=h_1gg_{\psi}h_1^{-1}$. Indeed, if $g_{\psi'} \notin A_i$, we are done. Otherwise, since $\Poly(\phi|_{\mathcal{F}})=\mathrm{Poly}(H|_{\mathcal{F}})$ and $\Poly(\phi) \neq \Poly(H)$, there exist $\psi \in H$ and $h \in A-A_i$ such that either $g_{\psi} \notin A_i$, or $\Psi_0(A)=A$ and $\Psi_0(A_i)=hA_ih^{-1}$. In the first case we are done. Otherwise, we have $\Psi_0 \circ \Psi_0'(g)=gg_{\psi}hg_{\psi'}h^{-1}$. Since $A_i$ is malnormal, we have $hg_{\psi}h^{-1} \notin A_i$ and $g_{\psi}hg_{\psi}h^{-1} \notin A_i$. The claim follows. Thus, for every $n \in \NN^*$ and $\psi \in H$ as in the claim, we have $g_{\psi^n} \notin A_i$ and there exists $h_n \in A$ such that $\Psi^n(A_i)=A_i$ and $\Psi^n(g)=h_ngg_{\psi^n}h_n^{-1}$. 

\medskip

\noindent{\bf Claim~2. } For every $n \in \NN^*$ and every $a,h \in F_{\n}$, if $$h \in (a\Psi^n(A_i \ast \left\langle g \right\rangle)a^{-1}) \cap (A_i \ast \left\langle g \right\rangle),$$ then $h$ is contained in a conjugate of $A_i$.

\medskip

\dem Let $n \in \NN^*$, let $a \in F_{\n}$ and let $h \in (a\Psi^n(A_i \ast \left\langle g \right\rangle)a^{-1}) \cap (A_i \ast \left\langle g \right\rangle)$. Suppose towards a contradiction that $h$ is not contained in a conjugate of $A_i$. First note that $ah_n gg_{\psi^n}h_n^{-1}a^{-1} \notin a A_i \ast \left\langle g \right\rangle a^{-1}$. Indeed, since $F_{\n}=A \ast \left\langle g \right\rangle$, the elements $h_n gg_{\psi^n}h_n^{-1}$ can be written uniquely as a reduced product of elements in $A$ and elements in $\left\langle g \right\rangle$. Since $h_n, g_{\psi^n} \in A$, if we have $ah_n gg_{\psi^n}h_n^{-1}a^{-1} \in  a A_i \ast \left\langle g \right\rangle a^{-1}$, then $h_n \in A_i$ and $g_{\psi^n}h_n^{-1} \in A_i$. Therefore, $g_{\psi^n} \in A_i$, a contradiction. Thus, we have $ah_n gg_{\psi^n}h_n^{-1}a^{-1} \notin  a A_i \ast \left\langle g \right\rangle a^{-1}$.

We claim that there exist a subgroup $B'$ of $F_{\n}$ containing $A_i \ast \left\langle g \right\rangle$, $aA_i \ast \left\langle g \right\rangle a^{-1}$ and $ah_ngg_{\psi^{n}}h_n^{-1}a^{-1}$, and an $\RR$-tree $T'$ equipped with a minimal, isometric action of $B'$ with trivial arc stabilizers and such that $A_i \ast \left\langle g \right\rangle$ fixes a point $x_1'$ in $T'$, $aA_i \ast \left\langle g \right\rangle a^{-1}$ fixes a point $x_2'$ in $T'$ and $ah_ngg_{\psi^n}h_n^{-1}a^{-1}$ either fixes a point in $T'$ distinct from $x_1'$ or $x_2'$ or is loxodromic. Indeed, by~\cite[Lemma~1.2]{Levitt09}, there exists a nontrivial $\RR$-tree $T$ equipped with a minimal, isometric action of $F_{\n}$ with trivial arc stabilizers and such that every polynomial subgroup of $\phi$ fixes a point in $T$. In particular, $A_i \ast \left\langle g \right\rangle$ fixes a point $x_1$ in $T$ and $aA_i \ast \left\langle g \right\rangle a^{-1}$ fixes a point $x_2$ in $T$. If $ah_ngg_{\psi^n}h_n^{-1}a^{-1}$ either fixes a point in $T$ distinct from $x_1$ or $x_2$ or is loxodromic, we may take $T=T'$. Otherwise $x_1=x_2$, $ah_ngg_{\psi^n}h_n^{-1}a^{-1} \in \Stab(x_1)$ and an induction on the rank of $\Stab(x_1)$ (which is less than $\n$ by~\cite{GabLev95} and invariant by a power of $\phi$) allows us to conclude since $ah_ngg_{\psi^n}h_n^{-1}a^{-1} \notin  a A_i \ast \left\langle g \right\rangle a^{-1}$. 

Suppose first that $ah_ngg_{\psi^n}h_n^{-1}a^{-1}$ fixes $x_2'$. Then $x_1' \neq x_2'$. Since $a\Psi^n(A_i \ast \left\langle g \right\rangle)a^{-1}$ fixes $x_2'$, $h$ fixes both $x_1'$ and $x_2'$. This contradicts the fact that $T'$ has trivial arc stabilizers.

Suppose now that $ah_ngg_{\psi^n}h_n^{-1}a^{-1}$ does not fix $x_2'$. We claim that $h$ is loxodromic in $T'$. Indeed, note that, since $h \in a\Psi^n(A_i \ast \left\langle g \right\rangle)a^{-1}$, $h$ can be written as a product of elements of $aA_ia^{-1}$ and powers of $ah_ngg_{\psi^n}h_n^{-1}a^{-1}$. Since $h$ is not contained in a conjugate of $A_i$, we may suppose that:

\noindent{\hypertarget{1}{$(i)$} } the word $h$ contains at least one occurence of a nontrivial element in $aA_ia^{-1}$ and one occurrence of a nontrivial power of $ah_ngg_{\psi^n}h_n^{-1}a^{-1}$;

\noindent{\hypertarget{2}{$(ii)$} } the word $h$ is cyclically reduced when written in the generating set $$\{aA_ia^{-1},ah_ngg_{\psi^n}h_n^{-1}a^{-1}\}.$$ 

Suppose first that $ah_ngg_{\psi^n}h_n^{-1}a^{-1}$ fixes a point $x'$ in $T'$ (which is distinct from $x_2'$). Then the minimal tree $T_0'$ in $T'$ of the subgroup $B$ of $F_{\n}$ generated by $\Stab(x')$ and $\Stab(x_2')$ is simplicial and its vertex stabilizers are conjugates of $\Stab(x')$ and $\Stab(x_2')$. We conclude as in the proof of Claim~1 that $h$ is loxodromic. Suppose now that $ah_ngg_{\psi^n}h_n^{-1}a^{-1}$ is loxodromic and that its axis does not contain $x_2'$. Assertion~\hyperlink{2}{$(ii)$} implies that, if $h$ is not loxodromic, then it fixes $x_2'$. Then a ping pong argument shows, since $h$ satisfies Assertion~\hyperlink{1}{$(i)$}, that $h$ is loxodromic. Finally, suppose that $ah_ngg_{\psi^n}h_n^{-1}a^{-1}$ is loxodromic and that its axis contains $x_2'$. Assertion~\hyperlink{2}{$(ii)$} implies that, if $h$ is not loxodromic, then it fixes $x_2'$. Then the minimal tree $T_0'$ in $T'$ of the subgroup $B$ of $F_{\n}$ generated by $\Stab(x_2')$ and $ah_ngg_{\psi^n}h_n^{-1}a^{-1}$ is simplicial and its vertex stabilizers are conjugate of $\Stab(x_2')$ (it is an HNN extension). In particular, $h$ is loxodromic as it satisfies Assertion~\hyperlink{1}{$(i)$}. Thus $h$ is loxodromic in $T'$.  Hence $h$ cannot fix $x_1'$ and $h \notin A_i \ast \left\langle g \right\rangle$, a contradiction.
\hfill\qedsymbol  

\medskip

Claim~$2$ implies that, for every distinct $n,m \in \NN^*$ and every $x \in F_{\n}$, we have $$\Psi^n(A_i \ast \left\langle g \right\rangle) \cap x\Psi^m(A_i \ast \left\langle g \right\rangle)x^{-1}\subseteq \bigcup_{y \in F_{\n}} yA_iy^{-1}$$ and by~\cite[Fact~I.1.2]{HandelMosher20}, we have $$\partial^2\Psi^n(A_i \ast \left\langle g \right\rangle) \cap \partial^2\left(x\Psi^m(A_i \ast \left\langle g \right\rangle)x^{-1}\right)\subseteq \overline{\bigcup_{y \in F_{\n}} \partial^2 \left(yA_iy^{-1}\right)}. $$ The rest of the proof is then similar to the one of Case~1.
\hfill\qedsymbol

\begin{lem}\label{Lem preparation ping pong}
Let $\n \geq 3$, let $\mathcal{F}$ and $\mathcal{F}_1=\{[A_1],\ldots,[A_k]\}$ be two free factor systems of $F_{\n}$ with $\mathcal{F} \leq \mathcal{F}_1$ such that the extension $\mathcal{F} \leq \mathcal{F}_1$ is sporadic. Let $H$ be a subgroup of $\Out(F_{\n},\mathcal{F},\mathcal{F}_1) \cap \IA_{\n}(\ZZ/3\ZZ)$ such that $H$ is irreducible with respect to $\mathcal{F} \leq \mathcal{F}_1$. Suppose that there exists $\phi \in H$ such that $\Poly(\phi|_{\mathcal{F}})=\mathrm{Poly}(H|_{\mathcal{F}})$. Suppose that \mbox{$\Poly(\phi|_{\mathcal{F}_1}) \neq \Poly(H|_{\mathcal{F}_1})$}. There exists $\psi \in H$ such that for every $i \in \{1,\ldots,k\}$, we have $\psi(K_{PG}([\Phi_i|_{A_i}])) \cap K_{PG}([\Phi_i|_{A_i}])=\varnothing$ and $$\Delta_+([A_i],\phi) \cap \psi(\Delta_-([A_i],\phi))=\Delta_-([A_i], \phi) \cap \psi(\Delta_+([A_i], \phi))=\varnothing.$$
\end{lem}

\dem The proof follows \cite[Lemma~5.1]{clay2019atoroidal}. Recall that, since the extension $\mathcal{F} \leq \mathcal{F}_1$ is sporadic, there exists $\ell \in \{1,\ldots,k\}$ such that, for every $j \in \{1,\ldots,k\}-\{\ell\}$, we have $[A_j] \in \mathcal{F}$. By Lemma~\ref{Lem Kpg moved by elements} applied to the image of $H$ in $\Out(A_{\ell})$ (which is contained in $\IA(A_{\ell},\ZZ/3\ZZ)$), there exists an infinite subset $X \subseteq H$ such that, for any distinct $h_1,h_2 \in X$, we have $$h_1(K_{PG}([\Phi_{\ell}|_{A_{\ell}}])) \cap h_2(K_{PG}([\Phi_{\ell}|_{A_{\ell}}]))=\varnothing.$$ We now prove that there exist $h_1,h_2 \in X$ such that $h_2^{-1}h_1$ satisfies the assertion of Lemma~\ref{Lem preparation ping pong}. Note that, for any distinct $h_1,h_2 \in X$, we have $h_2^{-1}h_1(K_{PG}([\Phi_{\ell}|_{A_{\ell}}])) \cap K_{PG}([\Phi_{\ell}|_{A_{\ell}}])=\varnothing$. Hence it suffices to find two distinct $h_1,h_2$ such that $\psi=h_2^{-1}h_1$ satisfies the second assertion of Lemma~\ref{Lem preparation ping pong}. Let $j \in \{1,\ldots,k\}$ and let $[\mu]$ be an extremal point of $\Delta_+([A_j],\phi)$ or $\Delta_-([A_j],\phi)$. By \cite[Lemma~4.13]{Guerch2021NorthSouth}, the support $\Supp(\mu)$ contains the support of finitely many projective currents $[\mu_1],\ldots,[\mu_s] \in \PCurr(F_n,\mathcal{F}\wedge \mathcal{A}(\phi))$ such that, for every $t \in \{1,\ldots,s\}$, the support of $\mu_t$ is uniquely ergodic. Let $E_{\mu}=\{[\mu_1],\ldots,[\mu_s]\}$. Let $E_{\phi}=\bigcup E_{\mu}$, where the union is taken over all $j \in \{1,\ldots,k\}$ and extremal points of $\Delta_+(A_j,\phi)$ and $\Delta_-(A_j,\phi)$. The set $E_{\phi}$ is finite by \cite[Lemma~4.7]{Guerch2021NorthSouth}.

Since the set $E_{\phi}$ is finite, up to taking an infinite subset of $X$, we may suppose that, for every $s \in E_{\phi}$, either $h_1s=h_2s$ for every $h_1,h_2 \in X$ or for every distinct $h_1,h_2 \in X$, we have $h_1s \neq h_2s$. Let $E_1 \subseteq E_{\phi}$ be the subset for which the first alternative occurs and let $E_{\infty}=E_{\phi}-E_1$. 

Let $h_1 \in X$ and, for every $s \in E_{\infty}$, let 
$$X_s=\{h \in X \;|\; h_1s=hs' \text{ for some }s' \in 
E_{\infty}\}.$$  Note that $X_s$ is a finite set. Let $h_2 \in X-\bigcup_{s \in E_{\infty}} X_s$. For every $s,s' \in E_{\infty}$, we have $h_1s \neq h_2 s'$. If there exists $s' \in E_1$ such that $h_1s=h_2s'$, then $s=h_1^{-1}h_2s'=s'$, contradicting the fact that $s \in E_{\infty}$. Thus, for every $s \in E_{\infty}$, we have $h_2^{-1}h_1s \notin E_{\phi}$ and for every $s \in E_1$, we have $h_2^{-1}h_1s=s$. Let $\psi=h_2^{-1}h_1$. Then, for every $s \in E_{\phi}$, either $\psi(s)=s$ or $\psi(s) \notin E_{\phi}$.

Let $j \in \{1,\ldots,k\}$, let $[\mu] \in \Delta_-([A_j], \phi)$ and suppose for a contradiction that we have $\psi([\mu]) \in \Delta_+([A_j],\phi)$. There exist extremal measures $\mu_1^-,\ldots,\mu_m^-$ of $\Delta_-([A_j],\phi)$ and $\lambda_1,\ldots,\lambda_m \in \RR_+$ such that $\mu=\sum_{i=1}^m \lambda_i\mu_i^-$. Similarly, there exist extremal measures $\mu_1^+,\ldots,\mu_n^+$ of $\Delta_+([A_j],\phi)$ and $\alpha_1,\ldots,\alpha_n \in \RR_+$ such that $\psi(\mu)=\sum_{i=1}^n \alpha_i\mu_i^+$.

Thus, we have $$\sum_{i=1}^m \lambda_i\psi(\mu_i^-)=\psi(\mu)=\sum_{i=1}^n \alpha_i\mu_i^+.$$ In particular, we have $$\bigcup_{i=1}^m\Supp(\psi(\mu_i^-))=\bigcup_{i=1}^n\Supp(\mu_i^+).$$ Let $\Lambda \subseteq \Supp(\mu_1^-)$ be the uniquely ergodic support of a current in $E_{\phi}$. Let $\Psi$ be a representative of $\psi$ and let $\partial^2\Psi$ be the homeomorphism of $\partial^2 F_{\n}$ induced by $\Psi$. Since uniquely ergodic laminations are minimal, there exists $i \in \{1,\ldots,n\}$ such that we have $\partial^2\Psi(\Lambda) \subseteq \Supp(\mu_i^+)$. Thus, we have $\psi([\mu_1^-|_{\Lambda}])=[\mu_i^+|_{\Lambda}]$. This contradicts the fact that $[\mu_1^-|_{\Lambda}]$ and $[\mu_i^+|_{\Lambda}]$ are distinct elements of $E_{\phi}$ since $\Delta_+([A_j],\phi) \cap \Delta_-([A_j],\phi)=\varnothing$.
\hfill\qedsymbol

\begin{prop}\label{Prop ping pong}
Let $\n \geq 3$, let $\mathcal{F}$ and $\mathcal{F}_1=\{[A_1],\ldots,[A_k]\}$ be two free factor systems of $F_{\n}$ with $\mathcal{F} \leq \mathcal{F}_1$ such that the extension $\mathcal{F} \leq \mathcal{F}_1$ is sporadic. Let $H$ be a subgroup of $\IA_{\n}(\ZZ/3\ZZ) \cap \Out(F_{\n},\mathcal{F},\mathcal{F}_1)$ such that $H$ is irreducible with respect to $\mathcal{F} \leq \mathcal{F}_1$. Suppose that there exists $\phi \in H$ such that $\Poly(\phi|_{\mathcal{F}})=\mathrm{Poly}(H|_{\mathcal{F}})$. Suppose that $\Poly(\phi|_{\mathcal{F}_1}) \neq \Poly(H|_{\mathcal{F}_1})$. There exist $\psi \in H$ and a constant $M>0$ such that, for all $m,n \geq M$, if $\theta=\psi\phi\psi^{-1}$, we have $\Poly(\theta^m\phi^n|_{\mathcal{F}_1})=\Poly(H|_{\mathcal{F}_1})$.
\end{prop}

\dem The proof follows~\cite[Proposition~5.2]{clay2019atoroidal}. Let $\psi \in H$ be an element given by Lemma~\ref{Lem preparation ping pong} and let $\theta=\psi\phi\psi^{-1}$. For every $i \in \{1,\ldots,k\}$, let $\Theta_i$ be a representative of $\theta$ such that $\Theta_i(A_i)=A_i$ and $\Phi_i$ be a representative of $\phi$ such that $\Phi_i(A_i)=A_i$. Note that, since for every $i \in \{1,\ldots,k\}$, $[\Phi_i|_{A_i}]$ is almost atoroidal relative to $\mathcal{F}$, so is $[\Theta_i|_{A_i}]$. Moreover, for every $i \in \{1,\ldots,k\}$, we have $K_{PG}([\Theta_i|_{A_i}])=[\Psi_i|_{A_i}](K_{PG}([\Phi_i|_{A_i}]))$. Let $i \in \{1,\ldots,k\}$. Let $A_i \wedge \mathcal{F}$ be the free factor system of $A_i$ induced by $\mathcal{F}$: it is the free factor system of $A_i$ consisting in the intersection of $A_i$ with every subgroup $A$ of $F_{\n}$ such that $[A] \in \mathcal{F}$. It is well-defined by for instance~\cite[Theorem~3.14]{ScoWal79}. 

\medskip

\noindent{\bf Claim. } We have $\widehat{\Delta}_+([A_i],\phi) \cap \psi(\widehat{\Delta}_-([A_i],\phi))=\varnothing$ and $\widehat{\Delta}_-([A_i],\phi) \cap \psi(\widehat{\Delta}_+([A_i],\phi))=\varnothing$. 

\medskip

\dem We prove the first equality, the other one being similar. By Lemma~\ref{Lem preparation ping pong}, we have $\Delta_+([A_i],\phi) \cap \psi(\Delta_-([A_i],\phi))=\varnothing$ and $\psi(K_{PG}([\Phi_i|_{A_i}])) \cap K_{PG}([\Phi_i|_{A_i}])=\varnothing$. Let $[\mu] \in \widehat{\Delta}_+([A_i],\phi) \cap \psi(\widehat{\Delta}_-([A_i],\phi))$. By definition, there exist $[\mu_1] \in \Delta_+([A_i],\phi)$, $[\nu_1] \in K_{PG}([\Phi_i|_{A_i}])$, $t \in [0,1]$, and $[\mu_2] \in \psi(\Delta_-([A_i],\phi)$, $[\nu_2] \in \psi(K_{PG}([\Phi_i|_{A_i}]))$, $s \in [0,1]$ such that $$[\mu]=[t\mu_1+(1-t)\nu_1]=[s\mu_2+(1-s)\nu_2].$$ Note that $$\partial^2 (\mathcal{F} \wedge \{[A_i]\}) \cap \partial^2\mathcal{A}(\phi) \cap \partial^2(A_i, \mathcal{F} \wedge \{[A_i]\}\wedge \mathcal{A}(\phi))=\varnothing.$$ Moreover, since $\Poly(\phi|_{\mathcal{F}})=\Poly(H|_{\mathcal{F}})$, we have $\Poly(\theta|_{\mathcal{F}})=\Poly(H|_\mathcal{F})$. Therefore, we see that $\mathcal{F}\wedge\mathcal{A}(\phi)=\mathcal{F}\wedge \psi(\mathcal{A}(\phi))$. Thus, we have $$\partial^2 (\mathcal{F} \wedge \{[A_i]\}) \cap \psi(\partial^2\mathcal{A}(\phi)) \cap \partial^2(A_i, \mathcal{F} \wedge \{[A_i]\}\wedge \mathcal{A}(\phi))=\varnothing.$$ Recall that, by Proposition~\ref{Prop Existence and properties of Delta}, the supports of the currents in $$\Delta_+([A_i],\phi) \cup \psi(\Delta_-([A_i],\phi))$$ are contained in $\partial^2(\mathcal{F}\wedge \{[A_i]\})$. Moreover, by definition, the supports of currents in $K_{PG}([\Phi_i|_{A_i}])$ are contained in $\partial^2\mathcal{A}(\phi) \cap \partial^2(A_i, \mathcal{F} \wedge \{[A_i]\}\wedge \mathcal{A}(\phi))$ and the supports of currents in $\psi(K_{PG}([\Phi_i|_{A_i}]))$ are contained in $\psi(\partial^2\mathcal{A}(\phi)) \cap \partial^2(A_i, \mathcal{F} \wedge \{[A_i]\}\wedge \mathcal{A}(\phi))$. Thus, we have $$\mu_1(\partial^2\mathcal{A}(\phi) \cap \partial^2(A_i, \mathcal{F} \wedge \{[A_i]\}\wedge \mathcal{A}(\phi)))=\mu_2(\partial^2\mathcal{A}(\phi) \cap \partial^2(A_i, \mathcal{F} \wedge \{[A_i]\}\wedge \mathcal{A}(\phi)))=0.$$ Hence the support of $\nu_1$ is contained in the support of $\nu_2$. By definition of $\psi(K_{PG}([\Phi_i|_{A_i}]))$, this implies that $$\nu_1 \in K_{PG}([\Phi_i|_{A_i}]) \cap \psi(K_{PG}([\Phi_i|_{A_i}]))=\varnothing.$$ Thus, we necessarily have $t=1$. Similarly, we have $s=1$. This implies that $[\mu_1]=[\mu_2]$ and that $\Delta_+([A_i],\phi) \cap \psi(\Delta_-([A_i],\phi))\neq \varnothing,$ a contradiction. 
\hfill\qedsymbol

\medskip

By the claim, there exist subsets $U,V, \widehat{U}, \widehat{V}$ of $\PCurr(A_i,(A_i \wedge \mathcal{F}) \wedge \mathcal{A}(\phi))$ such that:

\medskip

\noindent{$(1)$ } $\Delta_+([A_i],\phi) \subseteq U$, $\widehat{\Delta}_+([A_i],\phi) \subseteq \widehat{U}$, $\Delta_-([A_i],\phi) \subseteq V$, $\widehat{\Delta}_-([A_i],\phi) \subseteq \widehat{V}$;

\noindent{$(2)$ } $U \subseteq \widehat{U}$, $V \subseteq \widehat{V}$;

\noindent{$(3)$ } $\widehat{U} \cap \psi (\widehat{V})=\varnothing$ and $\widehat{V} \cap \psi (\widehat{U})=\varnothing$.

\medskip

Let $\mathfrak{B}$ and $C>0$ be respectively the basis of $F_{\n}$ and the constant given by Proposition~\ref{Prop moving exponential1}~$(1)$. Let $M_0(\phi)$ (resp. $M_0(\theta^{-1})$) be the constant associated with $\phi$, $U$ and $\widehat{V}$ (resp $\theta^{-1}$, $\psi(V)$ and $\psi(\widehat{U})$) given by Theorem~\ref{Theo North-South dynamics almost atoroidal}. Let $M_1(\phi)$ and $L_1(\phi)$, (resp. $M_1(\theta)$ and $L_1(\theta)$) be the constants associated with $[\Phi_i|_{A_i}]$ and $\widehat{V}$ (resp. $[\Theta_i|_{A_i}]$ and $\psi(\widehat{V})$) given by Proposition~\ref{Prop moving exponential}. Similarly, let $M_1(\phi^{-1})$ and $L_1(\phi^{-1})$, (resp. $M_1(\theta^{-1})$ and $L_1(\theta^{-1})$) be the constants associated with $[\Phi_i|_{A_i}^{-1}]$ and $\widehat{U}$ (resp. $[\Theta_i|_{A_i}^{-1}]$ and $\psi(\widehat{U})$) given by Proposition~\ref{Prop moving exponential}. Let $$M=\max\{M_0(\phi),M_0(\theta^{-1}),M_1(\phi),M_1(\theta),M_1(\phi^{-1}),M_1(\theta^{-1})\}$$ and let $$L=\min\{L_1(\phi),L_1(\theta),L_1(\phi^{-1}),L_1(\theta^{-1})\} >0.$$

Let $M'$ be such that $3^{M'}L^2>1$. Let $m,n \geq M+M'$ and let $\mu \in \Curr(A_i,A_i \wedge \mathcal{F} \wedge \mathcal{A}(\phi))$ be a nonzero current. Suppose first that $[\mu] \notin \widehat{V}$. Then by Theorem~\ref{Theo North-South dynamics almost atoroidal}, we have $\phi^n(\mu) \in U$. By Proposition~\ref{Prop moving exponential}, we have $\lVert \phi^n(\mu) \rVert_{\mathcal{F}} \geq 3^{n-M}L\lVert \mu \rVert_{\mathcal{F}}$. Since $U \cap \psi(\widehat{V})=\varnothing$, by Proposition~\ref{Prop moving exponential}, we have

$$\lVert \theta^m\phi^n(\mu) \rVert_{\mathcal{F}} \geq 3^{m-M}L\lVert \phi^n(\mu) \rVert_{\mathcal{F}} \geq 3^{m+n-2M}L^2\lVert \mu \rVert_{\mathcal{F}}.$$ Note that, since $\widehat{V} \cap \psi(\widehat{U})=\varnothing$, we have $\theta^m\phi^n(\mu) \notin \widehat{V}$. Therefore, we can apply the same arguments replacing $\mu$ by $\theta^m\phi^n(\mu)$ and an inductive argument shows that, for every $n' \in \NN^*$, we have $$\lVert (\theta^m\phi^n)^{n'}(\mu) \rVert_{\mathcal{F}} \geq 3^{n'(m+n-2M-M')}(3^{M'}L^2)^{n'}\lVert \mu \rVert_{\mathcal{F}}.$$ Therefore, if $\mu$ is the current associated with a nonperipheral element $g \in A_i$ with $[\mu] \notin \widehat{V}$, for every $n' \geq 1$, by Proposition~\ref{Prop moving exponential1}~$(1)$ we have $$\ell_{\mathfrak{B}}((\theta^m\phi^n)^{n'}([g])) \geq 3^{n'(m+n-2M-M')}(3^{M'}L^2)^{n'}C\lVert \mu \rVert_{\mathcal{F}} \geq 3^{n'(m+n-2M-M')}C.$$ Hence we have $g \notin \Poly([\Theta_i^m\Phi_i^n|_{A_i}])$. Suppose now that $[\mu] \in \widehat{V}$. Therefore, we have $[\mu] \notin \psi(\widehat{U})$. By Theorem~\ref{Theo North-South dynamics almost atoroidal}, we have $\theta^{-m}([\mu]) \in \psi(V)$. By Proposition~\ref{Prop moving exponential}, we have $\lVert \theta^{-m}(\mu) \rVert_{\mathcal{F}} \geq 3^{m-M}L\lVert \mu \rVert_{\mathcal{F}}$. Moreover, since $\psi(V) \cap \widehat{U}=\varnothing$, we have $\theta^{-m}([\mu]) \notin \widehat{U}$ and $$\lVert \phi^{-n}\theta^{-m}(\mu) \rVert_{\mathcal{F}} \geq 3^{n-M}L\lVert \theta^{-m}(\mu) \rVert_{\mathcal{F}} \geq 3^{n+m-2M-M'}(3^{M'}L^2)\lVert \mu \rVert_{\mathcal{F}}.$$ Note that, since $\widehat{V} \cap \psi(\widehat{U})=\varnothing$, we have $\phi^{-n}\theta^{-m}(\mu) \notin \psi(\widehat{U})$. Therefore, we can apply the same arguments replacing $\mu$ by $\phi^{-n}\theta^{-m}(\mu)$ and an inductive argument shows that, for every $n' \in \NN^*$, we have $$\lVert (\phi^{-n}\theta^{-m})^{n'}(\mu) \rVert_{\mathcal{F}} \geq 3^{n'(m+n-2M-M')}(3^{M'}L^2)^{2n'}\lVert \mu \rVert_{\mathcal{F}}.$$ Therefore, if $\mu$ is the current associated with a nonperipheral element $g \in A_i$ with $[\mu] \in \widehat{V}$, for every $n' \geq 1$, we have $$\ell_{\mathfrak{B}}((\phi^{-n}\theta^{-m})^{n'}([g])) \geq 3^{n'(m+n-2M-M')}(3^{M'}L^2)^{n'}C\lVert \mu \rVert_{\mathcal{F}} \geq 3^{n'(m+n-2M)}C.$$ Hence we have $g \notin \Poly([\Phi_i^{-n}\Theta_i^{-m}|_{A_i}])=\Poly([\Theta_i^m\Phi_i^n|_{A_i}])$. Therefore, $\theta^m\phi^n|_{\mathcal{F}_1}$ is expanding relative to $\mathcal{F} \wedge \mathcal{A}(\phi)$. Hence if $g \in \Poly(\theta^m\phi^n|_{\mathcal{F}_1})$, there exists a subgroup $A$ of $F_{\n}$ such that $g \in A$ and $[A] \in \mathcal{F} \wedge \mathcal{A}(\phi)$. Note that, if $g \in F_{\n}$ is such that there exists a subgroup $A$ of $F_{\n}$ such that $g \in A$ and $[A] \in \mathcal{F} \wedge \mathcal{A}(\phi)$, then $g \in \Poly(\phi|_{\mathcal{F}})=\Poly(H|_\mathcal{F})$. Thus, we have $\Poly(\theta^m\phi^n|_{\mathcal{F}_1})=\Poly(H|_{\mathcal{F}_1})$. This concludes the proof.
\hfill\qedsymbol

\begin{prop}\label{Prop pushing sporadic extensions}
Let $\n \geq 3$ and let $H$ be a subgroup of $\IA_{\n}(\ZZ/3\ZZ)$. Let $$\varnothing=\mathcal{F}_0 < \mathcal{F}_1 <\ldots < \mathcal{F}_k=\{[F_{\n}]\} $$ be a maximal $H$-invariant sequence of free factor systems. Let $2 \leq i \leq k$. Suppose that $\mathcal{F}_{i-1} \leq \mathcal{F}_i$ is sporadic. Suppose that there exists $\phi \in H$ such that

\noindent{$(a)$ } $\Poly(H|_{\mathcal{F}_{i-1}})=\Poly(\phi|_{\mathcal{F}_{i-1}})$;

\noindent{$(b)$ } for every $j \in \{1,\ldots,k\}$, if the extension $\mathcal{F}_{j-1} \leq \mathcal{F}_j$ is nonsporadic, then $\phi|_{\mathcal{F}_j}$ is fully irreducible relative to $\mathcal{F}_{j-1}$ and if $H|_{\mathcal{F}_{j}}$ is atoroidal relative to $\mathcal{F}_{j-1}$, so is $\phi|_{\mathcal{F}_j}$.

Then there exists $\widehat{\phi} \in H$ such that:

\noindent{$(1)$ } $\Poly(H|_{\mathcal{F}_{i}})=\Poly(\widehat{\phi}|_{\mathcal{F}_{i}})$;

\noindent{$(2)$ } for every $j \in \{1,\ldots,k\}$, if the extension $\mathcal{F}_{j-1} \leq \mathcal{F}_j$ is nonsporadic, then $\widehat{\phi}|_{\mathcal{F}_j}$ is fully irreducible relative to $\mathcal{F}_{j-1}$ and if $H|_{\mathcal{F}_{j}}$ is atoroidal relative to $\mathcal{F}_{j-1}$, so is $\widehat{\phi}|_{\mathcal{F}_j}$.
\end{prop}

\dem The proof follows~\cite[Proposition~5.3]{clay2019atoroidal}. If $\Poly(H|_{\mathcal{F}_{i}})=\Poly(\phi|_{\mathcal{F}_{i}})$, we may take $\widehat{\phi}=\phi$. Otherwise, by Proposition~\ref{Prop ping pong}, there exists $\psi \in H$ and a constant $M>0$ such that, for every $m,n \geq M$, if $\theta=\psi\phi\psi^{-1}$, we have $\Poly(\theta^m\phi^n|_{\mathcal{F}_i})=\Poly(H|_{\mathcal{F}_i})$. Therefore, for every $m,n \geq M$, the element $\widehat{\phi}=\theta^m\phi^n$ satisfies~$(1)$. It remains to show that there exist $m,n \geq M$ such that $\theta^m\phi^n$ satisfies $(2)$. Let $$S=\{j\;|\; \text{the extension } \mathcal{F}_{j-1} \leq \mathcal{F}_j \text{ is nonsporadic}\}$$ and let $j \in S$. 

Let $B_j$ be a subgroup of $F_{\n}$ given by the claim in the proof of Theorem~\ref{Theo nonsporadic extension}. Let $A_{j,1},\ldots,A_{j,s}$ be the subgroups of $B_j$ with pairwise disjoint conjugacy classes such that $\mathcal{A}_{j-1}=\{[A_{j,1}],\ldots,[A_{j,s}]\} \subseteq \mathcal{F}_{j-1}$ and $s$ is maximal for this property. By Hypothesis~$(b)$, the outer automorphism $[\Phi_j|_{B_j}] \in \Out(B_j,\mathcal{A}_{j-1})$ is fully irreducible relative to $\mathcal{A}_{j-1}$. Note that $\mathcal{A}_{j-1}$ is a nonsporadic free factor system of $B_j$ by the claim and since the extension $\mathcal{F}_{j-1} \leq \mathcal{F}_j$ is nonsporadic. Let $S_1$ be the subset of $S$ consisting in every $j \in S$ such that $H|_{\mathcal{F}_j}$ is atoroidal relative to $\mathcal{F}_{j-1}$, and let $S_2=S-S_1$. By~\cite[Theorem~4.1,4.2]{Guirardelhorbez19} (see also~\cite{Mann2014,Mann2014Thesis,Horbez17cyclic,gupta18}), for every $j \in S_1$ (resp. $j \in S_2$) there exists a Gromov-hyperbolic space $X_j$ (the \emph{$\mathcal{Z}$-factor complex of $B_j$ relative to $\mathcal{A}_{j-1}$} when $j \in S_1$ and the \emph{free factor complex of $B_j$ relative to $\mathcal{A}_{j-1}$} otherwise) on which $\Out(B_j,\mathcal{A}_{j-1})$ acts by isometries and such that $\phi_0 \in \Out(B_j,\mathcal{A}_{j-1})$ is a loxodromic element if and only if $\phi_0$ is fully irreducible atoroidal relative to $\mathcal{A}_{j-1}$ (resp. fully irreducible relative to $\mathcal{A}_{j-1}$). In particular, since $H$ preserves $\mathcal{F}_{j-1} <\mathcal{F}_j$, and hence $\psi$ preserves $[B_j]$, the elements $[\Phi_j|_{B_j}]$ and $[\Theta_j|_{B_j}]$ are loxodromic elements of $X_j$. 

Recall that two loxodromic isometries of a Gromov-hyperbolic space $X$ are \emph{independent} if their fixed point sets in $\partial_{\infty} X$ are disjoint and are \emph{dependent} otherwise. Let $I \subseteq S$ be the subset of indices where for every $j \in I$, the elements $[\Phi_j|_{B_j}]$ and $[\Theta_j|_{B_j}]$ are independent and let $D=S-I$. By standard ping pong arguments (see for instance~\cite[Proposition~4.2, Theorem~3.1]{ClayUya2018}), there exist constants $m,n_0 \geq M$ such that for every $n \geq n_0$, the element $[\Theta_j^m\Phi_j^n|_{B_j}]$ acts loxodromically on $X_j$. By~\cite[Proposition~3.4]{ClayUya2018}, there exists $n \geq n_0$ such that, for every $j \in D$, the element $[\Theta_j^m\Phi_j^n|_{B_j}]$ acts loxodromically on $X_j$. Hence for every $j \in S_1$, the element $\theta^m\phi^n|_{\mathcal{F}_j}$ is fully irreducible atoroidal relative to $\mathcal{F}_{j-1}$ and for every $j \in S_2$, the element $\theta^m\phi^n|_{\mathcal{F}_j}$ is fully irreducible relative to $\mathcal{F}_{j-1}$. This concludes the proof.
\hfill\qedsymbol

\section{Proof of the main result and applications}

We are now ready to complete the proof of our main theorem.

\begin{theo}\label{Theo poly growth}
Let $\n \geq 3$ and let $H$ be a subgroup of $\Out(F_{\n})$. There exists $\phi \in H$ such that $\Poly(\phi)=\Poly(H)$.
\end{theo}

\dem Since $\IA_{\n}(\ZZ/3\ZZ)$ is a finite index subgroup of $\Out(F_{\n})$ and since for every $\psi \in H$ and every $n \in \NN^*$, we have $\Poly(\psi^k)=\Poly(\psi)$, we see that $\Poly(H)=\Poly(H \cap \IA_{\n}(\ZZ/3\ZZ))$. Hence we may suppose that $H$ is a subgroup of $\IA_{\n}(\ZZ/3\ZZ)$. Let $$\varnothing=\mathcal{F}_0 < \mathcal{F}_1 <\ldots < \mathcal{F}_k=\{[F_{\n}]\} $$ be a maximal $H$-invariant sequence of free factor systems. By Theorem~\ref{Theo nonsporadic extension}, there exists $\phi \in H$ such that for every $j \in \{1,\ldots,k\}$ such that the extension $\mathcal{F}_{j-1} \leq \mathcal{F}_j$ is nonsporadic, the element $\phi|_{\mathcal{F}_j}$ is fully irreducible relative to $\mathcal{F}_{j-1}$ and if $H|_{\mathcal{F}_j}$ is atoroidal relative to $\mathcal{F}_{j-1}$, so is $\phi|_{\mathcal{F}_{j-1}}$.

We now prove by induction on $i \in \{0,\ldots,k\}$ that for every $i \in \{0,\ldots,k\}$, there exists $\phi_i \in H$ such that

\medskip

\noindent{$(a)$ } $\Poly(\phi_i|_{\mathcal{F}_i})=\Poly(H|_{\mathcal{F}_i})$;

\medskip

\noindent{$(b)$ } for every $j \in \{1,\ldots,k\}$ such that the extension $\mathcal{F}_{j-1} \leq \mathcal{F}_j$ is nonsporadic, the element $\phi_i|_{\mathcal{F}_j}$ is fully irreducible relative to $\mathcal{F}_{j-1}$ and if $H|_{\mathcal{F}_j}$ is atoroidal relative to $\mathcal{F}_{j-1}$, so is $\phi_i|_{\mathcal{F}_{j-1}}$. 

\medskip

For the base case $i=0$, we set $\phi_0=\phi$. Let $i \in \{1,\ldots,k\}$ and suppose that $\phi_{i-1} \in H$ has been constructed. We distinguish between two cases, according to the nature of the extension $\mathcal{F}_{i-1} \leq \mathcal{F}_i$. Suppose first that the extension $\mathcal{F}_{i-1} \leq \mathcal{F}_i$ is nonsporadic. We set $\phi_i=\phi_{i-1}$. We claim that $\phi_i$ satisfies the hypotheses. Indeed, it clearly satisfies $(b)$. For $(a)$, since $\Poly(\phi_{i-1}|_{\mathcal{F}_{i-1}})=\Poly(H|_{\mathcal{F}_{i-1}})$, it suffices to show that for every element $g \in F_{\n}$ which is $\mathcal{F}_i$-peripheral but $\mathcal{F}_{i-1}$-nonperipheral, if $g \in \Poly(\phi_i|_{\mathcal{F}_i})$, then $g \in \Poly(H|_{\mathcal{F}_i})$. Note that, if $\phi_i|_{\mathcal{F}_i}$ is atoroidal relative to $\mathcal{F}_{i-1}$, by Proposition~\ref{Prop properties fully irreducible}~$(1)$, we have $\Poly(\phi_i|_{\mathcal{F}_i})=\Poly(\phi_i|_{\mathcal{F}_{i-1}})$. Hence we have $\Poly(H|_{\mathcal{F}_i})=\Poly(\phi_i|_{\mathcal{F}_i})$. So we may suppose that $\phi_i|_{\mathcal{F}_i}$ is not atoroidal relative to $\mathcal{F}_{i-1}$. 

Let $g \in \Poly(\phi_i|_{\mathcal{F}_i})$ be an element which is $\mathcal{F}_i$-peripheral but $\mathcal{F}_{i-1}$-nonperipheral. By Proposition~\ref{Prop properties fully irreducible}~$(1)$, there exists at most one (up to taking inverse) $h \in F_{\n}$ such that $g \in \left\langle h \right\rangle$ and $[h]$ is fixed by $\phi_i$. By Proposition~\ref{Prop properties fully irreducible}~$(2)(b)$, the conjugacy class of $[h]$ is fixed by $H$. Hence the conjugacy class of $[g]$ is fixed by $H$ and $g \in \Poly(H|_{\mathcal{F}_i})$.

Suppose now that the extension $\mathcal{F}_{i-1} \leq \mathcal{F}_i$ is sporadic. If $\Poly(\phi_{i-1}|_{\mathcal{F}_i})=\Poly(H|_{\mathcal{F}_i})$, we set $\phi_i=\phi_{i-1}$. Then $\phi_i$ satisfies $(a)$ and $(b)$. Suppose that $\Poly(\phi_{i-1}|_{\mathcal{F}_i})\neq \Poly(H|_{\mathcal{F}_i})$. By Proposition~\ref{Prop pushing sporadic extensions}, there exists $\widehat{\phi}_{i-1} \in H$ such that $\widehat{\phi}_{i-1}$ satisfies $(a)$ and $(b)$. Then we set $\phi_i=\widehat{\phi}_{i-1}$. This completes the induction argument. In particular, we have $\Poly(\phi_m)=\Poly(H)$. This concludes the proof.
\hfill\qedsymbol

\bigskip

We now give some applications of Theorem~\ref{Theo poly growth}. The first one is a straightforward consequence using the fact that for every $\phi \in \Out(F_{\n})$, there exists a natural malnormal subgroup system associated with $\Poly(\phi)$.

\begin{coro}\label{Coro malnormal subgroup system}
Let $\n \geq 3$ and let $H$ be a subgroup of $\Out(F_n)$ such that $\Poly(H) \neq \{1\}$. There exist nontrivial maximal subgroups $A_1,\ldots,A_k$ of $F_{\n}$ such that $$\Poly(H)=\bigcup_{i=1}^k\bigcup_{g \in F_{\n}} gA_ig^{-1} $$ and $\mathcal{A}=\{[A_1],\ldots,[A_k]\}$ is a malnormal subgroup system. 
\hfill\qedsymbol
\end{coro}

If $H$ is a subgroup of $\Out(F_{\n})$ is such that $\Poly(H)\neq \{1\}$, we denote  by $\mathcal{A}(H)$ the malnormal subgroup system given by Corollary~\ref{Coro malnormal subgroup system}. If $\Poly(H)=\{1\}$, we set $\mathcal{A}(H)=\varnothing$.

The following result is a generalization of \cite[Theorem~A]{clay2019atoroidal} regarding fixed conjugacy classes. If $\phi\in \IA_{\n}(\ZZ/3\ZZ)$, we denote by $\mathrm{Fix}(\phi)$ the set of conjugacy classes of $F_{\n}$ fixed by $\phi$. Note that, if $g \in F_{\n}$ is such that $[g] \in \mathrm{Fix}(\phi)$, then $g \in \Poly(\phi)$. Moreover, by \cite[Lemma~1.5]{Levitt09}, if $\Poly(\phi) \neq \{1\}$, the set $\mathrm{Fix}(\phi)$ is nonempty. If $P$ is a subgroup of $F_{\n}$, we denote by $\Out(F_{\n},P^{(t)})$ the subgroup of $\Out(F_{\n})$ consisting in every element $\phi \in \Out(F_{\n})$ such that there exists $\Phi \in \phi$ such that $\Phi(P)=P$ and $\Phi|_{P}=\id_P$.

\begin{coro}\label{Coro alternative single elements}
Let $\n \geq 3$ and let $g_1,\ldots,g_k$ be nontrivial root-free elements of $F_{\n}$. Let $H$ be a subgroup of $\IA_{\n}(\ZZ/3\ZZ)$ such that, for every $i \in \{1,\ldots,k\}$, every element of $H$ fixes the conjugacy class of $g_i$. Then one of the following (mutually exclusive) statements holds.

\medskip

\noindent{$(1)$ } There exists $g_{k+1} \in F_{\n}$ such that $[g_{k+1}] \notin \{[g_1],\ldots,[g_k]\}$ and whose conjugacy class is fixed by every element of $H$.

\medskip

\noindent{$(2)$ } There exists $\phi \in H$ such that $\mathrm{Fix}(\phi)=\{[\left\langle g_1 \right\rangle],\ldots,[\left\langle g_k \right\rangle]\}$.

Moreover, if $(1)$ holds, either there exist $\ell \geq k+1$ and $g_{k+1},\ldots,g_{\ell} \in F_{\n}$ such that $$\mathrm{Fix}(H)=\mathcal{A}(H)=\{[\left\langle g_1 \right\rangle],\ldots,[\left\langle g_{\ell} \right\rangle]\}$$ or $H$ virtually fixes the conjugacy class of a nonabelian free subgroup of $F_{\n}$ of rank $2$.
\end{coro}

\dem First assume that $H$ is finitely generated. Suppose that $(2)$ does not hold. In particular, by Theorem~\ref{Theo poly growth}, we see that $\mathcal{A}(H) \neq \{[\left\langle g_1 \right\rangle],\ldots,[\left\langle g_k \right\rangle]\}$. Let $\mathcal{A}(H)=\{[P_1],\ldots,[P_{\ell}]\}$, where for every $i \in \{1,\ldots,\ell\}$, $P_i$ is a malnormal subgroup of $F_{\n}$. Note that, for every $i \in \{1,\ldots,\ell\}$, since $P_i$ is malnormal, we have a natural homomorphism $H \to \Out(P_i)$ whose image, denoted by $H|_{P_i}$, is contained in the set of polynomially growing outer automorphisms of $P_i$. Since $H$ is finitely generated, up to taking a finite index subgroup of $H$, we can apply the Kolchin theorem for $\Out(F_{\n})$ (see~\cite[Theorem~1.1]{BesFeiHan05}): there exists a $H|_{P_i}$-invariant sequence of free factor systems of $P_i$ $$\varnothing=\mathcal{F}_0^{(i)} < \mathcal{F}_1^{(i)} <\ldots < \mathcal{F}_{k_i}^{(i)}=\{[P_i]\} $$ such that, for every $j \in \{1,\ldots,k_i\}$, the extension $\mathcal{F}_{j-1}^{(i)} \leq \mathcal{F}_j^{(i)}$ is sporadic. 

Note that, since $\mathcal{A}(H) \neq \{[\left\langle g_1 \right\rangle],\ldots,[\left\langle g_k \right\rangle]\}$, either $\ell>k$ or there exists $i \in \{1,\ldots,\ell\}$ such that the rank of $P_i$ is at least equal to $2$. Suppose that $\ell >k$. Let $i \in \{1,\ldots,\ell\}$. Since for every $j \in \{1,\ldots,k_i\}$, the extension $\mathcal{F}_{j-1}^{(i)} \leq \mathcal{F}_j^{(i)}$ is sporadic, for every $i \in \{1,\ldots,\ell\}$, the free factor system $\mathcal{F}_1^{(i)}$ contains a unique element and the rank of the associated subgroup is $1$. Thus, the group $H$ fixes at least $\ell$ distinct conjugacy classes of elements of $F_{\n}$ and $(1)$ holds.

Otherwise, let $i \in \{1,\ldots,\ell\}$ be such that the rank of $P_i$ is at least equal to $2$. Since, for every $j \in \{1,\ldots,k_i\}$, the extension $\mathcal{F}_{j-1}^{(i)} \leq \mathcal{F}_j^{(i)}$ is sporadic we have $k_i \geq 2$. Moreover, there exists $j_0 \in \{1,\ldots,k_i\}$ and a subgroup $U_{j_0}$ of $P_i$ such that $[U_{j_0}] \in \mathcal{F}_{j_0}^{(i)}$ and one of the following holds:

\noindent{$(a)$ } there exist two subgroups $B_1$ and $B_2$ of $P_i$ such that $\mathrm{rank}(B_1)=\mathrm{rank}(B_2)=1$,  $[B_1],[B_2] \in \mathcal{F}_{j_0-1}$ and $U_{j_0}=B_1 \ast B_2$;

\noindent{$(b)$ } there exists a subgroup $B$ of $P_i$ such that $\mathrm{rank}(B)=1$,  $[B] \in \mathcal{F}_{j_0-1}$ and $U_{j_0}$ is an HNN extension of $B$ over the trivial group. 

If Case~$(a)$ occurs, then $H$ acts as the identity on $U_{j_0}$ since $\mathrm{rank}(U_{j_0})=2$ and since every element of $H$ fixes elementwise a set of conjugacy classes of generators of $U_{j_0}$ (recall that the abelianization homomorphism $F_2 \to \ZZ^2$ induces an isomorphism $\Out(F_2) \simeq \mathrm{GL}(2,\ZZ)$). Hence Assertion~$(1)$ holds. 

If Case~$(b)$ occurs, let $b$ be a generator of $B$ and let $t \in U_{j_0}$ be such that $U_{j_0}=\left\langle b \right\rangle \ast \left\langle t \right\rangle$. Then, since $H \subseteq \IA_{\n}(\ZZ/3\ZZ)$, for every element $\psi$ of $H$, there exist $\Psi \in \psi$ and $k \in \ZZ$ such that $\psi(b)=b$ and $\psi(t)=tb^k$. In particular, for every $\psi \in H$, the automorphism $\Psi$ fixes the group generated by $b$ and $tbt^{-1}$ and $(1)$ holds. 

The moreover part follows since either for every $i \in \{1,\ldots,\ell\}$, the group $P_i$ has rank $1$ or there exists $i \in \{1,\ldots,k\}$ such that the rank of $P_i$ is at least equal to $2$. In the first case, since $H \subseteq \IA_{\n}(\ZZ/3\ZZ)$, for every $i \in `\{1,\ldots,\ell\}$, the conjugacy class $[P_i]$ is fixed by $H$. In the later case, the subgroup $H$ fixes the conjugacy class of a nonabelian subgroup of rank $2$ as explained above. This concludes the proof when $H$ is finitely generated.

Suppose now that $H$ is not finitely generated and let $(H_m)_{m \in \NN}$ be an increasing sequence of finitely generated subgroups of $H$ such that $H=\bigcup_{m \in \NN} H_m$. For every $m \in \NN$, we have $H_m \subseteq \Out(F_{\n},\mathrm{Fix}(H_m)^{(t)})$ and for every $m_1,m_2 \in \NN$ such that $m_1 \leq m_2$, we have $\mathrm{Fix}(H_{m_2}) \subseteq \mathrm{Fix}(H_{m_1})$. By \cite[Theorem~1.5]{GuirardelLevitt2015McCool}, there exists $N \in \NN$ such that, for every $m \geq N$, we have $\Out(F_{\n},\mathrm{Fix}(H_m)^{(t)})=\Out(F_{\n},\mathrm{Fix}(H_N)^{(t)})$. In particular, we have $\mathrm{Fix}(H_N)=\mathrm{Fix}(H)$. The result now follows from the finitely generated case.
\hfill\qedsymbol

\bigskip

The following result might be folklore, but we did not find a precise statement in the literature. If $S$ is a compact, connected surface, we denote by $\Mod(S)$ the group of homotopy classes of homeomorphisms that preserve the boundary of $S$.

\begin{coro}\label{Coro pseudo Anosov subgroup Mod}
Let $\n \geq 3$ and let $H$ be a subgroup of $\IA_{\n}(\ZZ/3\ZZ)$. The following assertions are equivalent:

\medskip

\noindent{$(1)$ } $\mathcal{A}(H)=\{[\left\langle g \right\rangle]\}$, where $g$ is an element of $F_{\n}$ not contained in a proper free factor of $F_{\n}$;

\medskip

\noindent{$(2)$ } there exists a connected, compact surface $S$ with exactly one boundary component and an identification of $\pi_1(S)$ with $F_{\n}$ such that $H$ is identified with a subgroup of $\Mod(S)$ and $H$ contains a pseudo-Anosov element.
\end{coro}

\dem Suppose that $(2)$ holds. Let $\phi \in H$ be identified with a pseudo-Anosov element of $S$. In particular, $\phi$ is a fully irreducible element of $\Out(F_{\n})$. By Proposition~\ref{Prop properties fully irreducible}~$(1)$ with $\mathcal{F}=\varnothing$, the element $\phi$ fixes exactly one (up to taking inverse) conjugacy class $[g]$ of a root-free element $g$ of $F_{\n}$. Since $\phi$ fixes the conjugacy class of the element of $F_{\n}$ identified with the boundary component of $S$, the conjugacy class $[g]$ is identified with the conjugacy class in $\pi_1(S)$ of the element associated with the homotopy class of the boundary component of $S$. Hence $g$ is not contained in any proper free factor of $F_{\n}$. Moreover, since $H$ is identified with a subgroup of $\Mod(S)$, every element of $H$ fixes $[g]$. Hence we have $\mathcal{A}(H)=\{[\left\langle g \right\rangle]\}$. 

Suppose now that $(1)$ holds. Let $\phi \in H$ be an element given by Theorem~\ref{Theo poly growth}. Then $\mathcal{A}(\phi)=\mathcal{A}(H)=\{[\left\langle g \right\rangle]\}$. In particular, since $H \subseteq \IA_{\n}(\ZZ/3\ZZ)$, the conjugacy class of $g$ is fixed by every element of $H$. Let $f \colon G \to G$ be a CT map representing a power of $\phi$ (see the definition in \cite[Definition~4.7]{FeiHan06}). 

\medskip

\noindent{\bf Claim. } The graph $G$ consists in a single stratum and this stratum is an EG stratum. 

\medskip

\dem Let $H_r$ be the highest stratum in $G$. We first prove that $H_r$ is an EG stratum. Indeed, $H_r$ is either a zero stratum, an EG stratum or a NEG stratum. The stratum $H_r$ cannot be a zero stratum by \cite[Definition~4.7~$(6)$]{FeiHan06}. Moreover, $H_r$ cannot be a NEG stratum as otherwise by \cite[Proposition~4.1]{clay2019atoroidal}, the element $g$ would be a basis element of $F_{\n}$, contradicting the fact that $g$ is not contained in any proper free factor of $F_{\n}$. Hence $H_r$ is an EG stratum. Since $g$ is not contained in any proper free factor of $F_{\n}$, the reduced circuit $\gamma_g$ in $G$ representing the conjugacy class of $g$ has height $r$ and is fixed by $f$. By \cite[Fact~I.2.3]{HandelMosher20}, the stratum $H_r$ is a geometric stratum. By~\cite[Proposition~I.2.18]{HandelMosher20}, the element $\phi$ fixes elementwise a finite set $\mathcal{C}=\{[g],[c_1],\ldots,[c_k]\}$ of conjugacy classes of elements of $F_{\n}$. Moreover, by~\cite[Proposition~I.2.18~$(5)$]{HandelMosher20}, by the definition of a geometric stratum in~\cite{HandelMosher20} and the fact that $G$ is connected, we have $\mathcal{C}=\{[g]\}$ if and only if $G_{r-1}$ is reduced to a point, that is, if and only if $G$ consists in the single stratum $H_r$.
\hfill\qedsymbol

\medskip

By the claim and~\cite[Fact~I.2.3]{HandelMosher20}, the outer automorphism $\phi$ is \emph{geometric}: there exist a connected, compact surface $S$ with exactly one boundary component and an identification of $\pi_1(S)$ with $F_{\n}$ such that $\phi$ is identified with a pseudo-Anosov element of $\Mod(S)$. Moreover, the conjugacy class $[g]$ is identified with the conjugacy class in $\pi_1(S)$ of the element associated with the homotopy class of the boundary component of $S$. Since $[g]$ is fixed by every element of $H$, by the Dehn-Nielsen-Baer theorem (see for instance~\cite[Theorem~8.8]{FarMar12} for the orientable case and \cite[Section~3]{Fujiwara2002} for the nonorientable case), the group $H$ is identified with a subgroup of $\Mod(S)$.
\hfill\qedsymbol

\bigskip

We now give a method to compute $\mathcal{A}(H)$ for some subgroups $H$ of $\Out(F_{\n})$. Let $P$ be a subgroup of $F_{\n}$. Let $F$ be the minimal free factor of $F_{\n}$ which contains $P$. Then $F$ is one-ended relative to $P$. Let $T$ be the JSJ tree of $F$ relative to $P$ over cyclic subgroups given by \cite[Theorem~9.14]{guirardel2016jsj}. Let $v$ be a vertex of $T$. Let $G_v$ be the stabilizer of $v$ in $F$. Let $\mathrm{Inc}_v$ be the finite set of all conjugacy classes of groups associated with edges in $T$ which are incident to $v$. Following the terminology of \cite{guirardel2016jsj}, either $v$ is a \emph{rigid} vertex or $v$ is \emph{flexible}. When $G_v$ is cyclic, we use the convention that $v$ is rigid. If $v$ is flexible, by~\cite[Theorem~9.14~$(2)$]{guirardel2016jsj}, there exists a compact connected hyperbolic surface $S_v$ such that $\pi_1(S_v)$ is isomorphic to $G_v$ and, for every subgroup $G_e$ of $F$ such that $[G_e] \in \mathrm{Inc}_v$, the group $G_e$ is conjugate to a subgroup of $\pi_1(S_v)$ associated with a boundary connected component of $S_v$. Since the JSJ tree constructed by Guirardel and Levitt is a \emph{tree of cylinders}, if $v$ is a flexible vertex of $T$, the fundamental group of every boundary component $c$ of $S_v$ fixes at most one edge $e_c$ adjacent to $v$ and the stabilizer of the endpoint of $e_c$ distinct from $v$ is cyclic and included in the group generated by the homotopy class of $c$. For every flexible vertex $v$ of $T$, let $\mathrm{BC}_v$ be the finite set of conjugacy classes of subgroups of $F_{\n}$ generated by the homotopy classes of the boundary components of $S_v$ which do not fix an edge in $T$. Let $V_f$ be the set of flexible vertices of $T$. Let $T'$ be the tree obtained from $T$ by collapsing every edge of $T$ which is not adjacent to a flexible vertex. For every vertex $C$ of $F \backslash (T'-V_f)$, let $G_C$ be the associated vertex stabilizer. Let $$\mathcal{A}_P=\{[G_C]\}_{C \in V(F \backslash (T'-V_f))} \cup \bigcup_{v \in V_f} (\mathrm{Inc}_v\cup \mathrm{BC}_v),$$ which is a finite set of conjugacy classes of finitely generated subgroups of $F_{\n}$. Note that, by~\cite[Proposition~2.1]{Levitt09}, if $v \in VT$ and $\phi \in \Out(G_v)$, then $\phi$ extends to an outer automorphism $\widehat{\phi}$ of $\Out(F_{\n})$ which preserves $T$ and $T'$.

\begin{lem}\label{Lem exp growth extended pseudo Anosov}
Let $v \in V_f$ and let $\phi'$ be the outer automorphism class of $\Out(G_v)$ associated with a pseudo-Anosov element of $S_v$. Let $T_v'$ be the tree obtained from $T$ by collapsing every edge of $T$ which is not contained in the orbit of an edge adjacent to $v$. Then $\widehat{\phi'}$ preserves $T_v'$. Moreover, if $g \in F$ is loxodromic in $T_v'$, then $g$ has exponential growth under iteration of $\widehat{\phi'}$.
\end{lem}

\dem The fact that $\widehat{\phi'}$ preserves $T_v'$ follows from the fact that $\widehat{\phi'}$ preserves $T$ and the fact that $\widehat{\phi'}$ acts as the identity on the graph associated with $F_{\n} \backslash T$. In order to prove the second part of Lemma~\ref{Lem exp growth extended pseudo Anosov}, we first construct an $\RR$-tree $T_v$ with an equivariant map $T_v \to T_v'$. Let $T_0$ be the dilating arational $G_v$-tree associated with $\phi'$ and let $\lambda >1$ be the stretching factor of $\phi'$. There exists a homothety $H_0 \colon T_0 \to T_0$ whose stretching factor is equal to $\lambda$ and such that, for every $h \in G_v$ and every $x \in T_0$, we have $H_0(hx)=\phi'(h)H_0(x)$. The arational tree $T_0$ is such that every arc stabilizer is trivial and the only point stabilizers are cyclic and conjugate to the groups generated by the homotopy classes of boundary components of $S_v$. Since the edge stabilizers of $T_v'$ are precisely groups which are conjugates of groups generated by the homotopy classes of boundary components of $S_v$, one can replace the vertex $v$ in $T_v'$ by the tree $T_0$ and attach the edges to their corresponding point stabilizers. Extending this construction equivariantly, we obtain a tree $T_v$ with an equivariant map $T_v \to T_v'$. Moreover, since $\widehat{\phi'}$ preserves $T_v'$, the map $H_0$ extends to an equivariant map $H' \colon T_v \to T_v$ such that, for all $x,y \in T_v$, we have $d(H'(x),H'(y)) \leq \lambda d(x,y)$. We now follow the construction given by Gaboriau, Jaeger, Levitt and Lustig in \cite{GabJagLevLus98}. For every $k \in \NN^*$, let $d_k$ be the pseudo-distance in $T_v$ given by, for all $x,y \in T_v$: $$d_k(x,y)=\frac{d(H'^k(x),H'^k(y))}{\lambda^k},$$ and let $d_{\infty}$ be the limit of these pseudo-distances. Then $d_{\infty}$ induces a distance on the set $$T_v^{\infty}=T_v/\sim,$$ where $\sim$ is the equivalence relation generated by $x \sim y$ if and only if $d_{\infty}(x,y)=0$. Moreover, the metric space $T_v^{\infty}$ is a nontrivial $\RR$-tree equipped with a minimal, nontrivial action of $F_{\n}$ by isometries. Finally, $H'$ induces a homothety $H \colon T_v^{\infty} \to T_v^{\infty}$ with stretching factor equal to $\lambda$ and such that, for every $h \in F_{\n}$ and every $x \in T_v^{\infty}$, we have $H(hx)=\widehat{\phi'}(h)H(x)$. Note that, for every $g \in F_{\n}$ and every $n \in \NN^*$, the translation length of $\phi'^n(g)$ in $T_v^{\infty}$ is equal to $\lambda^n$ times the translation length of $g$ in $T_v^{\infty}$. Therefore, if $g$ has polynomial growth under iteration of $\phi'$, then $g$ must fix a point in $T_v^{\infty}$.

Let $g \in F$ be loxodromic in $T_v'$. By equivariance of the map $T_v \to T_v'$, the element $g$ is loxodromic in $T_v$. By the construction of $T_v'$, the axis of $g$ in $T_v$ contains a vertex in the orbit of $v$. Since the group generated by the homotopy class of every boundary component of $S_v$ fixes at most one edge in $T_v'$, if the axis of $g$ in $T_v'$ contains a vertex in the orbit of $v$, then the axis of $g$ in $T_v$ contains a nondegenerate arc $[x,y]$ in a copy of $T_0$ in $T_v$. Since $H_0$ is a homothety of $T_0$ of stretching factor equal to $\lambda$, the homothety $H$ restricts to a homothety of stretching factor $\lambda$ in the copy of $T_0$ in $T_v$. Thus we have $d_{\infty}(x,y)=d(x,y)>0$. Hence the characteristic set of $g$ in $T_v^{\infty}$, which is the projection of the characteristic set of $g$ in $T_v$, contains a nondegenerate arc, that is, $g$ is loxodromic in $T_v^{\infty}$. Hence $g$ has exponential growth under iteration of $\phi'$.  
\hfill\qedsymbol

\begin{prop}\label{Coro subgroup Fn are poly}
Let $\n \geq 3$ and let $P$ be a finitely generated subgroup of $F_{\n}$. Suppose that $\Out(F_{\n},P^{(t)})$ is infinite. Then either $\mathcal{A}(\Out(F_{\n},P^{(t)}))=\{[F_{\n}]\}$ or we have $\mathcal{A}(\Out(F_{\n},P^{(t)}))=\mathcal{A}_P$. Moreover, in the second case, there exists $\phi \in \Out(F_{\n},P^{(t)})$ such that $\mathcal{A}(\phi)=\mathcal{A}_P$.
\end{prop}

\dem The moreover part follows from the first part of Proposition~\ref{Coro subgroup Fn are poly} using Theorem~\ref{Theo poly growth}, so we focus on the first part. Let $H=\Out(F_{\n},P^{(t)}) \cap \IA_{\n}(\ZZ/3\ZZ)$. Then $\Poly(H)=\Poly(\Out(F_{\n},P^{(t)}))$.

Note that $H$ preserves the conjugacy class of $F$ and we have an induced homomorphism $\Lambda \colon H \to \Out(F,P^{(t)})$. By~\cite[Theorem~9.14]{guirardel2016jsj}, the group $\Lambda(H)$ fixes the $F$-equivariant homeomorphism class $\mathcal{T}$ of the above JSJ tree $T$. Moreover, up to taking a finite index subgroup of $H$, we may suppose that the group $H$ fixes the conjugacy class of every vertex group of $T$ and that $\Lambda(H)$ acts as the identity on the graph associated with $F \backslash T$. Since edge stabilizers of $T$ are cyclic, $H$ fixes the conjugacy class of the generator of every edge group. In particular, we have $\bigcup_{v \in V_f}\mathrm{Inc}_v \leq \mathcal{A}(H)$.  Moreover, up to taking a finite index subgroup of $H$, for every flexible vertex $v$ of $T$, the group $H$ fixes the conjugacy classes of subgroups of $F_{\n}$ generated by the homotopy classes of the boundary components of $S_v$. Thus, we have 
\begin{equation}\label{Equation inc included in Ah}
\bigcup_{v \in V_f}\mathrm{Inc}_v \cup \mathrm{BC}_v\leq \mathcal{A}(H). 
\end{equation}

\medskip

\noindent{\bf Claim~1. } Let $C'$ be a connected component of $F\backslash (T-V_f)$ and let $C$ be a connected subgraph of $C'$ which contains at least one vertex of the graph associated with $F\backslash T$. Let $F\backslash (T_C-V_f)$ be the graph of groups obtained from $F\backslash (T-V_f)$ by collapsing the edges of $F \backslash T$ contained in $C$ to a vertex $c$ and let $G_c$ be the corresponding vertex group. We have $G_c \subseteq \Poly(H)$. 

\medskip

\dem An \emph{interior edge of $C$} is an edge of the graph associated with $F \backslash (T-V_f)$ entirely contained in $C$. We remark that the statement of the claim is made in such a way that we are able to apply an induction argument on the number $m$ of interior edges of $C$. If $C$ does not contain an interior edge, then $C$ contains at most one vertex $v$ of $F \backslash T$ (recall that $C$ is connected). Moreover, $v$ is a rigid vertex. Since $v$ is rigid, by \cite[Theorem~3.9]{GuirardelLevitt2015}, the group $H$ has trivial image in $\Out(G_v)$. Hence the statement is true when $C$ has no interior edge. Suppose that the number of interior edges $m$ of $C$ is at least equal to $1$. Let $e$ be an interior edge of $C$. Suppose first that $\overline{C-e}$ has two connected components $A_1$ and $A_2$, where the closure is taken in $C$. For every $i \in \{1,2\}$, let $G_{a_i}$ be the subgroup of $F_{\n}$ corresponding to $A_i$ as in the statement of the claim. By induction, for every $i \in \{1,2\}$, we have $G_{a_i} \subseteq \Poly(H)$. Since $F$ is one-ended relative to $P$, edge stabilizers are nontrivial. Thus, we have $G_{a_1} \cap G_{a_2} \neq \{1\}$. Since $\mathcal{A}(H)$ is a malnormal subgroup system, there exists a subgroup $B$ of $F_{\n}$ such that $[B] \in \mathcal{A}(H)$ and $G_c=\left\langle G_{a_1},G_{a_2} \right\rangle \subseteq B$. Suppose now that $\overline{C-e}$ has one connected component $A$. Let $G_a$ be the subgroup of $F_{\n}$ corresponding to $A$ as in the statement of the claim. By induction, we have $G_{a} \subseteq \Poly(H)$. Moreover, there exists $t \in F$ such that $G_c=\left\langle G_a,t \right\rangle$. Note that $H$ preserves the conjugacy classes of $G_a$ and $G_c$ as every element of $\Lambda(H)$ acts as the identity on the graph associated with $F \backslash T$. Thus, every element $\psi$ of $H$ has a representative $\Psi$ such that $\Psi(G_a)=G_a$, every element of $G_a$ has polynomial growth under iteration of $\Psi$ and $\Psi$ sends $t$ to $ta_{\Psi}$ with $a_{\Psi} \in G_a$. Since $a_{\Psi}$ has polynomial growth under iteration of $\Psi$, there exist $s>0$, $n \in \NN$ and a free basis $\mathfrak{B}$ of $F_{\n}$ such that, for every $k \in \NN$, we have $$\ell_{\mathfrak{B}}(\Psi^k(a_{\Psi})) \leq s(k+1)^n.$$ Hence, for every $k \in \NN$, we have 
\begin{equation}\label{Equation polynomial growth HNN}
\ell_{\mathfrak{B}}(\Psi^k(t)) \leq 1+\sum_{i=0}^{k-1}\ell_{\mathfrak{B}}(\Psi^i(a_{\Psi}))\leq 1+ s(k+1)^{n+1}.
\end{equation} 
Therefore, $t$ has polynomial growth under iteration of $\Psi$ and every element of $G_c=\langle G_a,t\rangle$ has polynomial growth under iteration of $\Psi$. Thus, for every $\psi \in H$, we have $G_c \subseteq \Poly(\psi)$ and $G_c \subseteq \Poly(H)$. This proves the claim.
\hfill\qedsymbol

\medskip

By Claim~$1$ and Equation~\eqref{Equation inc included in Ah}, we have $\mathcal{A}_P \leq \mathcal{A}(H)$. We now prove that either $\mathcal{A}(H)=\{[F_{\n}]\}$ or $\mathcal{A}_P=\mathcal{A}(H)$. Let $K$ be a subgroup of $F_{\n}$ such that $F_{\n}=F \ast K$. Suppose first that the rank of $K$ is at most equal to $1$ and that the set of vertices $V_f$ is empty. Let $k$ be a (possibly trivial) generator of $K$. Recall the definition of $T'$ above Lemma~\ref{Lem exp growth extended pseudo Anosov}. Then $F \backslash T'$ is reduced to a vertex $v$. Therefore, we have $\mathcal{A}_P=\{[F]\}$. Moreover, since $H$ preserves the sporadic free factor system $\{[F]\}$, every element of $H$ has a representative which sends $F$ to $F$ and $k$ to $kg$ with $g \in F$. In particular, as in Equation~\eqref{Equation polynomial growth HNN}, we have $k \in \Poly(H)$, $F \ast K \subseteq \Poly(H)$ and $\mathcal{A}(H)=\{[F_{\n}]\}$. 

\medskip

\noindent{\bf Claim~2. } If either the rank of $K$ is at least equal to $2$ or $V_f$ is nonempty, then $\mathcal{A}(H) \leq \{[F]\}$.

\medskip

\dem We distinguish between two cases, according to the rank of $K$. When the rank of $K$ is equal to $0$, the proof is trivial.

\medskip

\noindent{\bf Case~1. } Suppose that the rank of $K$ is equal to $1$ and that $V_f$ is not empty. 

\medskip

Let $k$ be a generator of $K$. Let $v \in V_f$, let $\phi'$ be a pseudo-Anosov element of the surface $S_v$ associated with $\pi_1(G_v)$. As explained above Lemma~\ref{Lem exp growth extended pseudo Anosov}, the outer automorphism $\phi'$ induces an outer automorphism $\widehat{\phi'}$ of $F$. Let $g' \in G_v$ be such that $g'$ is not contained in the conjugacy class of the group generated by the homotopy class of  any boundary component of $S_{v}$. Let $\widehat{\Phi'}$ be a representative of $\widehat{\phi'}$. Let $\widehat{\Phi}$ be an automorphism of $F_{\n}$ which acts as $\widehat{\Phi'}$ on $F$ and sends $k$ to $kg'$, and let $\widehat{\phi}$ be the outer automorphism class of $\widehat{\Phi}$. Suppose that $\mathcal{A}(\widehat{\phi}) \nleq \{[F]\}$. By~\cite[Lemma~5.18~$(7)$]{Guerch2021NorthSouth} applied to $\mathcal{F}=\{[F]\}$ with the element $\widehat{\phi'} \in \Out(F_{\n},\mathcal{F})$ (recall that $F$ is a sporadic free factor of $F_{\n}$), there exists $g \in F_{\n}$ such that $F_{\n}=F \ast \left\langle g \right\rangle$ and either $\mathcal{A}(\widehat{\phi})=\mathcal{A}(\widehat{\phi'}) \cup \{[\langle g\rangle]\}$ or there exists a subgroup $A$ of $F$ such that $[A] \in \mathcal{A}(\widehat{\phi'})$ and $\mathcal{A}(\widehat{\phi})=(\mathcal{A}(\widehat{\phi'})-\{[A]\}) \cup \{[A \ast \langle g\rangle]\}$. In the first case, let $h \in F$ be nontrivial. Let $\Psi$ be the automorphism of $F_{\n}$ such that $\Psi(F)=F$, $\Psi|_F=\widehat{\Phi}'$ and $\Psi$ sends $g$ to $gh$ with $h \in F$ nontrivial and let $\psi$ be the outer automorphism class of $\Psi$. Note that $\psi \in \Out(F_{\n},P^{(t)})$. Then $\psi$ does not preserve the conjugacy class of $g$. Thus, we have $\mathcal{A}(H) \leq \mathcal{A}(\psi) \leq \{[F]\}$. 

Suppose that there exists a subgroup $A$ of $F$ such that $[A] \in \mathcal{A}(\widehat{\phi'})$ and $\mathcal{A}(\widehat{\phi})=(\mathcal{A}(\widehat{\phi'})-\{[A]\})\cup \{[A \ast \langle g\rangle ]\}$. Note that $\widehat{\phi}$ has a representative $\widehat{\Phi}_0$ such that $\widehat{\Phi}_0(F)=F$, $\widehat{\Phi}_0(A)=A$ and $\widehat{\Phi}_0(A \ast \langle g \rangle)=A \ast \langle g \rangle$. Then, up to composing $\widehat{\Phi}_0$ by an inner automorphism $\mathrm{ad}_{a_0}$ with $a_0 \in A$, we may suppose that $\widehat{\Phi}_0$ sends $g$ to $ga$ with $a \in A$. Moreover, since we have $\Poly(H) \subseteq \Poly(\widehat{\phi})$, if $h' \in F_{\n}$ is an $\{[F]\}$-nonperipheral element such that $h' \in \Poly(H)$, then $h'$ is contained in a conjugate of $A \ast \left\langle g \right\rangle$. Let $h \in F$ be an $\mathcal{A}(\widehat{\phi})$-nonperipheral element and let $\psi \in \Out(F_{\n},P^{(t)})$ be such that there exists $\Psi \in \psi$ with $\Psi(F)=F$, $\Psi|_F=\widehat{\Phi}_0|_F$ and $\Psi$ sends $g$ to $gh$. The element $h$ exists since $\phi'$ is a pseudo-Anosov of $S_v$, in particular, $\widehat{\phi'}$ is an exponentially growing outer automorphism of $F$. By Claim~2 in the proof of Lemma~\ref{Lem Kpg moved by elements}, for every $a' \in F_{\n}$, the intersection $(a'\Psi(A \ast \langle g \rangle)a'^{-1}) \cap (A \ast \langle g\rangle)$ is contained in conjugates of $A$. Note that $\Psi(A \ast \langle g \rangle)$ is the only (up to conjugacy) polynomial subgroup of $\psi\widehat{\phi}\psi^{-1}$ which contains $\{[F]\}$-nonperipheral element. Thus, every element of $A \ast \langle g\rangle$ which is not contained in a conjugate of $A$ has exponential growth under iteration of $\psi\widehat{\phi}\psi^{-1}$. In particular, we have $\mathcal{A}(H) \leq \{[F]\}$.

\medskip

\noindent{\bf Case~2. } Suppose that the rank of $K$ is at least equal to $2$.

\medskip

Note that we have $\Out(F_{\n},F^{(t)}) \subseteq \Out(F_{\n},P^{(t)})$ and that $F$ is a nonsporadic free factor of $F_{\n}$. By~\cite[Theorem~7.4]{Guirardelhorbez19}, since $\Out(F_{\n},F^{(t)})$ does not preserve the conjugacy class of any $\{[F]\}$-peripheral element of $F_{\n}$, the group $\Out(F_{\n},F^{(t)})$ contains a fully irreducible atoroidal element $\phi''$ of $F_{\n}$ relative to $\{[F]\}$. By Proposition~\ref{Prop properties fully irreducible}~$(1)$, there does not exist an $\{[F]\}$-nonperipheral element of $F_{\n}$ which has polynomial growth under iteration of $\phi''$. Thus we have $\mathcal{A}(\phi'') \leq \{[F]\}$. Thus, we have $\mathcal{A}(H) \leq \{[F]\}$. This proves Claim~2. 
\hfill\qedsymbol

\medskip

By Claim~$2$ and the paragraph above Claim~$2$, either $\mathcal{A}(H)=\{[F_{\n}]\}$ or $\mathcal{A}(H) \leq \{[F]\}$. We are thus left with the case $\mathcal{A}(H) \leq \{[F]\}$. In this case, we prove that $\mathcal{A}_P=\mathcal{A}(H)$. Since $\mathcal{A}_P \leq \mathcal{A}(H)$, it remains to prove that every $\mathcal{A}_P$-nonperipheral element of $F$ is $\mathcal{A}(H)$-nonperipheral. Let $g \in F$ be $\mathcal{A}_P$-nonperipheral. Recall that $T'$ is the tree obtained from $T$ by collapsing every edge of $T$ which is not adjacent to a flexible vertex. Note that, if a vertex $v$ of $T'$ is not the image of a flexible vertex of $T$, then its stabilizer is a conjugate of some $G_c$ with $C$ a connected component of $F \backslash (T-V_f)$. In particular, we have $[G_c] \in \mathcal{A}_P$. Suppose first that $g$ fixes a point in $T'$. Since $g$ is $\mathcal{A}_P$-nonperipheral, the element $g$ fixes a flexible vertex $v$ of $F_{\n}$ and is not conjugate to an element of $F_{\n}$ contained in the group generated by the homotopy class of a boundary component of $S_v$. Let $\phi' \in \Out(G_v)$ be the outer automorphism associated with a pseudo-Anosov element of $S_v$. Then $g$ has exponential growth under iteration of $\widehat{\phi'} \in \Out(F_{\n},P^{(t)})$. Thus, we have $g \notin \Poly(H)$. Suppose now that $g$ is loxodromic in $T'$. Then its axis contains the image of a flexible vertex $v \in V_f$. By Lemma~\ref{Lem exp growth extended pseudo Anosov}, $g$ has exponential growth under iteration of $\widehat{\phi'} \in \Out(F_{\n},P^{(t)})$. Therefore, every $\mathcal{A}_P$-nonperipheral element of $F$ is $\mathcal{A}(H)$-nonperipheral. Thus, we have $\mathcal{A}(H)=\mathcal{A}_P$.
\hfill\qedsymbol

\bibliographystyle{alphanum}
\bibliography{bibliographie}

\noindent \begin{tabular}{l}
Laboratoire de mathématique d'Orsay\\
UMR 8628 CNRS \\
Université Paris-Saclay\\
91405 ORSAY Cedex, FRANCE\\
{\it e-mail: yassine.guerch@universite-paris-saclay.fr}
\end{tabular}

\end{document}